\documentclass[journal,twoside]{IEEEtran}
\usepackage{graphicx}
\usepackage{textcomp}
\usepackage[usenames,dvipsnames,svgnames,table]{xcolor}
\usepackage[cmex10]{amsmath}
\usepackage{amssymb,amsfonts,mathrsfs,mathtools,amsthm}
\usepackage{epstopdf}
\usepackage{lscape}
\usepackage{CJK}
\usepackage{bm}
\usepackage{dsfont}
\usepackage{booktabs}
\usepackage{cite}
\usepackage{url}
\usepackage{subfigure}
\usepackage[colorlinks=true,allcolors=blue]{hyperref}
\usepackage{color}
\usepackage{algorithm,algorithmicx,algpseudocode}

\IEEEoverridecommandlockouts 

\newcommand{\tabincell}[2]{\begin{tabular}{@{}#1@{}}#2\end{tabular}} 

\definecolor{steelblue}{RGB}{70,130,180}
\definecolor{warmred}{RGB}{244, 104, 65}

\newtheorem{assumption}{Assumption}
\newtheorem{remark}{Remark}
\newtheorem{definition}{Definition}

\DeclareMathOperator*{\subjectto}{subject~to}

\begin{document}
\title{Optimal Pump Control for Water Distribution Networks via Data-based Distributional Robustness} 
\author{Yi Guo, \IEEEmembership{Member, IEEE}, 
Shen Wang,
Ahmad F. Taha, \IEEEmembership{Member, IEEE}
and Tyler H. Summers, \IEEEmembership{Member, IEEE}

\thanks{This material is based on works supported by the National Science Foundation under grants CMMI-1728605, CMMI-1728629, CMMI-2152450 and CMMI-2151392.}
\thanks{Y. Guo is with Department of Information Technology and Electrical Engineering at ETH Z\"{u}rich, Z\"{u}rich, 8092, Switzerland, email: yiguoy@ethz.ch.}
\thanks{S. Wang is with School of Cyberspace Security at Beijing University of Posts and Telecommunications, Beijing 100876, China, email: shen.wang@bupt.edu.cn.}
\thanks{A. F. Taha is with Department of Civil and Environmental Engineering at Vanderbilt University, Nashville, TN 37235, USA, email: ahmad.taha@vanderbilt.edu.}
\thanks{T. H. Summers is with the Department of Mechanical Engineering, The University of Texas at Dallas, Richardson,
TX 75080 USA, email: tyler.summers@utdallas.edu.}
}

\maketitle

\begin{abstract} 
In this paper, we propose a data-based methodology to solve a multi-period stochastic optimal water flow (OWF) problem for water distribution networks (WDNs). The framework explicitly considers the pump schedule and water network head level with limited information of demand forecast errors for an extended period simulation. The objective is to determine the optimal feedback decisions of network-connected components, such as nominal pump schedules and tank head levels and reserve policies, which specify device reactions to forecast errors for accommodation of fluctuating water demand. Instead of assuming the uncertainties across the water network are generated by a prescribed certain distribution, we consider \textit{ambiguity sets} of distributions centered at an empirical distribution, which is based directly on a finite training data set. We use a distance-based ambiguity set with the Wasserstein metric to quantify the distance between the real unknown data-generating distribution and the empirical distribution. This allows our multi-period OWF framework to trade off system performance and inherent sampling errors in the training dataset. Case studies on a three-tank water distribution network systematically illustrate the trade-off between pump operational cost, risks of constraint violation and out-of-sample performance. 
\end{abstract}


\section{Introduction}\label{sec:intro}
\IEEEPARstart{D}{ue} to a broad range of future energy and environmental issues \cite{water2035}, water distribution network operators are seeking improved strategies to deliver energy-efficient, reliable and high quality service to consumers \cite{ocampo2013application}. However, the increasing complexity (e.g., due to high dimensionality, nonlinearities, operation constraints and uncertainties) in municipal water supply network operation is challenging the current management and control strategies and may threaten the security of this vital infrastructure. Future urban water supply systems will require more sophisticated methods to function robustly and efficiently in the presence of this increasing complexity.

The flexibility of water flow manipulators (pumps and valves) in water networks has been utilized to optimize various objectives, including production and transportation costs, water quality, safe storage, smoothness of control actions etc. \cite{fooladivanda2018energy,d2015mathematical,jowitt1992optimal,yu1994optimized,verleye2013optimising,cohen2000optimal,montalvo2008particle,zessler1989optimal,wang2019geometric,wang2020receding}. However, most optimal water flow control methods use deterministic point forecasts of exogenous water demands, which neglects their inherent stochasticity. In practice, the variation of water demands in real water distribution networks is  high and difficult to predict \cite{goryashko2014robust}. Further, as complexity of network topology increases \cite{archibald2018review}, small perturbations can cause significant performance decrease and even infeasibility of optimal water flow problems \cite{goryashko2014robust}.

Recent research on optimal water network operation has been shifting from deterministic to stochastic models, since uncertainties (e.g., human usage, unexpected component failures, climate change) are increasingly key factors in many sectors of water resource management \cite{archibald2018review,brentan2018water,goryashko2014robust,watkins1997finding,sampathirao2018gpu,grosso2014economic,ocampo2013application,grosso2014chance,wang2016stochastic,castelletti2012stochastic,sankar2015optimal,wang2020probabilistic}. Most stochastic formulations assume that the uncertain water demands follow a prescribed distribution (e.g., Gaussian \cite{wang2016stochastic,grosso2014chance}), or enforce constraint for all possible water demand realization by assuming only knowledge of bounds on uncertainties \cite{goryashko2014robust,castelletti2012stochastic} and then utilize robust optimization. In addition, sampling-based stochastic optimization has also been applied to water flow manipulation problem \cite{sampathirao2018gpu} to quantify the probability of constraint violation based on an assumed data generating mechanism. However, the underlying assumptions in these approaches can be too strong or overly conservative, which can lead to underestimation or overestimation of the actual risks and therefore to ineffective operation. The methods based on chance-constraints effectively only measure the frequency of constraint violations not the severity, which can underestimate risk. The robust methods can enforce constraints for extreme and highly unlikely uncertainty realizations, effectively overestimating risk. Furthermore, some sampling-based methods are computational intensive due to their requirement of a large numbers of samples. In practice, forecasts of water demand are never perfect and their distributions must be estimated from finite data.


In this paper, we investigate a data-driven Model Predictive Control (MPC) approach to tackle a stochastic optimal water flow (OWF) problem for optimal pump schedule and head management in water distribution networks (WDNs). Many existing works, including \cite{camacho2013model,allgower2012nonlinear,rodriguez2021design,mayne2000constrained}, have successfully explored the stochastic MPC from the perspective of theoretical analysis and practical applications. We applied the proposed data-based distributionally robust stochastic MPC with the Wasserstein-based ambiguity set to address a water engineering problem. \color{black} Additionally, we also proposed a Laplacian-based water flow linearization associated with successive control algorithm, which is considered as a substantial contribution for dynamic control of water distribution networks from a control engineering perspectives. The proposed framework uses limited information of water demand forecasting errors from a finite training dataset to explicitly balance the trade-offs between performance and the risk of constraint violations in the presence of large water demand variations. A preliminary version of this work was presented in \cite{guo2020distributionally}; here we significantly extend the work in several directions. The main contributions are:

\begin{itemize}
\item[1)] We formulate a general multi-period distributionally robust optimal water flow problem for optimal pump schedule and head management. The distributionally robust OWF model predictive controller uses data-driven distributionally robust optimization \cite{mohajerin} to tractably obtain control decisions for network components at each stage. This allows the data-driven distributionally robust MPC OWF controller to update the water demand forecast with a finite time horizon and then re-compute the real-time optimal decision based on the latest and future forecasting information. In general, this OWF controller accepts the training data set from all forecasting frameworks and the decisions can be robust to various \textit{ambiguity sets (i.e., moment-based or metric-based)}. In this paper, we assume the unknown real data generating distribution is located in a metric-based ambiguity set, which is constructed by a Wasserstein ball with constant radius centered at an empirical distribution supported by the finite training dataset. In contrast to other stochastic OWF formulations, this approach makes the resulting control policies explicitly robust to the inherent sampling errors in the training dataset, which leads to superior out-of-sample performance. We can appropriately tune the size of the ambiguity set to avoid overly conservative decisions.

\item[2)]To handle computationally difficulties with the nonlinear/nonconvex water network hydraulics, we leverage a pertinent linear approximation of water network hydraulic coupling (i.e., flow-head coupling) to promote a computationally-efficient stochastic optimal water flow formulation for optimal pump control and nodal pressure management. In contrast with the literature, we further establish a generic matrix linearization in compact format between water flow and nodal head by re-defining a network Laplacian matrix based on linearization coefficients. This provides a unified framework that is applicable for approximation algorithms after linearization (i.e., successive linearization algorithms or piece-wise linearization algorithms). We empirically observe that the convergence of successive linearization algorithm provides an excellent approximation to the nonlinear water flow.

\item[3)] The effectiveness and flexibility of our proposed stochastic water flow formulation are demonstrated on a model of the Barcelona water distribution network. We illustrate the inherent trade-off between the system conservativeness and forecasting errors. The results can help the operators to explicitly prioritize the trade-off between the pump operational efficiency and the risk of tank head constraint violation and then develop the appropriate control strategies to balance their objectives and risk aversion.
\end{itemize}

The rest of paper is organized as follows: Section \ref{sec:mode and linearization} describes a generic model of water distribution networks and the successive linearization approach; Section \ref{sec:stochastic_OWF} presents the general formulation of proposed data-based multi-period distributionally robust stochastic optimal water flow problems. Section \ref{sec:data-driven_OWF} specifies the proposed stochastic OWF to a stochastic optimal pump schedule and head management. Section \ref{sec:case_studies} demonstrates the flexibility and effectiveness of the proposed methodologies via numerical experiments. Section \ref{sec:conclusion} concludes.\\
\textbf{Notation}: The inner product of two vectors $a,b \in \mathbf{R}^m$  is denoted by $\langle a, b \rangle := a^\intercal b$ and $(\cdot)^\intercal$ denotes the transpose of a  vector or matrix. The $N_s$-fold product of distribution $\mathds{P}$ on a set $\Xi$ is denoted by $\mathds{P}^{N_s}$, which represents a distribution on the Cartesian product space $\Xi^{N_s} = \Xi \times \ldots \times \Xi$. We use $N_s$ to represent the number of samples inside the training dataset $\hat{\Xi}$. Superscript `` $\hat{\cdot}$ " is reserved for the objects that depend on a training dataset $\hat{\Xi}^{N_s}$. The cardinality of set $\mathcal{J}$ is denoted by $|\mathcal{J}|$. The Kronecker product operator is defined as $\otimes$.

\section{Hydraulic Model and Leveraging Linear Approximation}\label{sec:mode and linearization}
In this section, we consider a water distribution network model associated with active and passive networked components and then we leverage a pertinent linear approximation, which leads to a novel network Laplacian-based matrix expression. This allows us to use successive linearization to approximate the original nonlinear hydraulic relationship for several topologies. WDNs control actions include speeds of pumps and settings of valves. In the rest of this section, we introduce the network and hydraulic modelling of networked components.

\subsection{Network Modelling}
We consider a water distribution network as a directed graph  
$\mathcal{G}(\mathcal{N},\mathcal{E})$ with a set $\mathcal{N}:=\{1,2,\ldots,N\}$ of vertices. These vertices include junctions, reservoirs and tanks that are collected in sets $\mathcal{J}$, $\mathcal{S}$ and $\mathcal{T}$ and $\mathcal{N} = \mathcal{J}\cup\mathcal{S}\cup\mathcal{T}$. Similarly, the set $\mathcal{E} \subseteq \mathcal{N}\times \mathcal{N}$ of all links including the sets of pipes, pumps and valves represented by $\mathcal{I}$, $\mathcal{M}$ and $\mathcal{V}$ so that $\mathcal{E} = \mathcal{I}\cup\mathcal{M}\cup\mathcal{V}$. Let $\mathcal{N}^{\textrm{in}}_i$ and $\mathcal{N}^{\textrm{out}}_i$ collect the supplying and carrying neighboring vertices of $i^\mathrm{th}$ node , respectively. We use $q_{ij} \in \mathbf{R}$ to denote the water flow through the link $(i,j) \in \mathcal{E}$ and $h_i \in \mathbf{R}_{+}$ denotes the head of node $i \in \mathcal{N}$. We assume each pipe has a prescribed flow direction, and that the sign of $q_{ij}$ indicates whether the actual flow direction matches the assumption ($q_{ij} \geq 0$) or goes in the opposite direction ($q_{ij} < 0$).

\subsubsection{Junctions}
The water demand is assumed to be a constant $d_i(t)$ in gallons per hour (GPM), which is applied for time interval $t$ at junction $i \in \mathcal{J}$. Mass conservation must be hold any time at $i^\mathrm{th}$ node  
\begin{equation}\label{eq:mass_balance}
    \sum_{j\in\mathcal{N}^{\textrm{in}}_i}q_{ji}(t) - \sum_{j\in\mathcal{N}^{\textrm{out}}_i}q_{ij}(t) = d_i(t), \quad \forall i \in \mathcal{J},
\end{equation}
where $\mathcal{N}^{\textrm{in}}_i$ and $\mathcal{N}^{\textrm{out}}_i$ are the sets of nodes supplying and carrying flow at $i^\mathrm{th}$ junction, respectively. If there is no water demand consumption for nodes $i \in \mathcal{N}\backslash \mathcal{N}_d$, we have $d_i(t) = 0$ for all time slots. Here we define an aggregated vector $d(t):=[d_1(t),\ldots,d_N(t)]^\intercal \in \mathbf{R}^{N}$.

\subsubsection{Reservoirs} The set $\mathcal{S}$ collects all reservoirs in a water distribution network. We assume that all reservoirs have infinite water resource supply and that the head of each reservoir is a constant, which can be treated as an operational constraint
\begin{equation}\label{eq:reservoir constraint}\nonumber
    h_i^{\textrm{R}} = h^{\textrm{elv}}_i, \quad \forall i\in \mathcal{S}, 
\end{equation}
where $h^{\textrm{elv}}_i$ represents the elevation for $i^\mathrm{th}$ reservoir.

\subsubsection{Tanks}
The head of tank at node $i \in\mathcal{T}$ at time $t$ is represented by $h_i^{\textrm{TK}}(t)$. The dynamics of these elements are given by the discrete-time difference equations
\begin{equation}\label{eq:tank_head_dynamics}
    h_{i}^{\textrm{TK}}(t+1)=h_{i}^{\textrm{TK}}(t)+\frac{\Delta t}{A_{i}^{\textrm{TK}}}\left(\sum_{j \in \mathcal{N}_{i}^{\mathrm{in}}} q_{ji}(t)-\sum_{j \in \mathcal{N}_{i}^{\mathrm{out}}} q_{ij}(t)\right),
\end{equation}
where $\Delta t$ is the duration of the time interval $(t, t+1]$. The cross-section area of tank is defined by $A_i^{\textrm{TK}}$.

\subsubsection{Pumps} The pumps provide head gain in the water distribution networks on the links $(i,j)\in\mathcal{M}$ connecting the suction $j^\mathrm{th}$ node  and the delivery $i^\mathrm{th}$ node. The head gain explicitly depends on the pump flow and pump property. Now we consider the variable speed pump (VSP) in the network and the head gain is given by
\begin{equation} \label{eq:pump_head_gain}
    h_{ij}^{\textrm{M}}(t) = h_i(t) - h_j(t) = \alpha_{ij}q_{ij}^{\textrm{M}}(t)^2 + \beta_{ij}q_{ij}^{\textrm{M}}(t) + \gamma_{ij},
\end{equation}
where coefficients $\alpha_{ij}$, $\beta_{ij}$ and $\gamma_{ij}$ are determined by the pump operation curve.

\subsubsection{Pipes} The head loss of pipe $(i,j)\in \mathcal{I}$ described via the empirical Chezy-Manning (C-M) is given as follows
\begin{equation} \label{eq:pipe_head_loss}
    h_{ij}^{\textrm{P}} = h_i(t) - h_j(t) = R_{\textrm{CM},ij}q_{ij}^{\textrm{P}}(t)^2,
\end{equation}
where the resistance coefficient is denoted by $R_{\textrm{CM},ij}\in\mathbf{R}_{++}$ and defined by \cite{rossman2000epanet}
\begin{equation}\label{eq:CM-coefficient}\nonumber
    R_{\textrm{CM},ij} = 4.66\frac{L_{\textrm{CM},ij}C^2_{\textrm{CM}}}{D_{\textrm{CM},ij}^{5.33}}.
\end{equation}
Note that $C_{\textrm{CM}}\in\mathbf{R}_{++}$ is the Manning roughness coefficient; $D_{\textrm{CM},ij}\in\mathbf{R}_{++}$ is the diameter of pipeline in feet; and $L_{\textrm{CM},ij}\in\mathbf{R}_{++}$ is the length of pipeline in feet.

\subsubsection{Pressure Reduce Valves} There are several types of controllable valves in a water distribution network, such as pressure reduce valves (PRVs), general purpose valves (GPVs) and flow control valves (FCVs), associated with different control variables: valve openness, pressure reduction and flow regulation. Here, we utilize PRVs to restrict the pressure to a certain difference $\phi_{ij}\in\mathbf{R}_{+}, (i,j)\in\mathcal{V}$ along a pipeline when the upstream pressure at $i^\mathrm{th}$ node is higher than the downstream $j^\mathrm{th}$ node 
\begin{equation}\label{eq:valve_head_reduce}
   \phi_{ij}(t) = h_i(t) - h_j(t),
\end{equation}
where $h_i(t)$ the head of upstream junction and $h_j(t)$ is the head of downstream junction and the variable $\phi_{ij}$ determines the energy conservation on pipeline $(i,j)$. Note that no reverse flow on PRVs is allowed and the water flow through PRVs, $q_{ij}, (i,j)\in\mathcal{V}$, is not determined by \eqref{eq:valve_head_reduce}, thereby depends on other network coupling constraints \eqref{eq:mass_balance}. The implementation of valve control actions depends on valve construction. We refer interested readers to \cite{fooladivanda2018energy,skworcow2014pump,prescott2003dynamic,al2006minimizing,ulanicki2008pressure} for more details. The deployment of PRVs in water distribution networks can promote the potential control availability. Here, we utilize a ``smart" PRV, whose pressure reduce setting can be optimized in the real time.

\subsubsection{Network Operational Constraints}
We specify several constraints on network states and inputs in our proposed stochastic OWF problem to satisfy the physical operation limitation of water distribution networks (i.e., limits of nodal heads, pipe flows and tank levels)
\begin{subequations}
\begin{equation}\label{eq:head_bounds}\nonumber
    h_i^{\textrm{min}} \leq h_i(t) \leq h_i^{\textrm{max}}, \quad \forall i \in \mathcal{N},
\end{equation}
\begin{equation}\label{eq:flow_bounds}\nonumber
    q_{ij}^{\textrm{min}} \leq q_{ij}(t) \leq q_{ij}^{\textrm{max}}, \quad  \forall (i,j) \in \mathcal{E},
\end{equation}
\end{subequations}
where $h_i^{\textrm{min}}$ and $h_i^{\textrm{max}}$ are the lower and upper heads on $i^\mathrm{th}$ node and $q_{ij}^{\textrm{min}}$ and $q_{ij}^{\textrm{max}}$ are minimum and maximum flows on link $(i,j)$. We introduce a binary parameter $z_{ij}(t)$ for link $(i,j)$ to indicate the ON/OFF status of the controllable devices (i.e., pumps, valves). Then the head coupling between two neighboring nodes can be modelled as follows
\begin{equation}\label{eq:general_energy_conservation}
\begin{aligned}
 -M\left(1 - z_{ij}(t)\right)~~~~~~~~~~~~~~~~~~~~~~~~~~~~~~~~~~~~~~ \\ \leq h_{ij}(t) - g(q_{ij}(t),\phi_{ij}(t)) \leq  M\left(1 - z_{ij}(t)\right),
 \end{aligned}
\end{equation}
where $h_{ij}(t):= h_i(t) - h_j(t)$, $g(\cdot)$ is a general expression of~\eqref{eq:pump_head_gain},~\eqref{eq:pipe_head_loss} and~\eqref{eq:valve_head_reduce} as functions of $q_{ij}$ and $\phi_{ij}$ and $M$ is a large positive constant. Note that when $z_{ij}(t) = 1$, the device on link $(i,j)$ is ON (i.e., $q_{ij}(t) \neq 0 $), then the energy conservation constraints \eqref{eq:pump_head_gain}, \eqref{eq:pipe_head_loss} and \eqref{eq:valve_head_reduce} hold on this link; otherwise $z_{ij}(t) = 0$ and the head at $i^\mathrm{th}$ node and $j^\mathrm{th}$ node are decoupled. For the links without a controllable device (e.g., pipes), we let $z_{ij}(t) = 1$ for all time intervals \footnote{Our water flow model for water distribution networks is given with binary variables in general. In the rest of this paper, we develop our distributionally robust MPC optimal water flow problem with the assumption that all binary variables are fixed, such that the ON/OFF status of the controllable devices are not optimized in our models.}.

\begin{table}[htbp!]
\begin{center}
\caption{Variable Notations}
\begin{tabular}{c|c}
\hline
\hline
Notation & Description\\
\hline
$h_i$  & Head at node $i$\\
\hline
$h_{i}^{\textrm{TK}}$/$h_{i}^{\textrm{R}}$  & Head at tank/reservoir $i$\\
\hline
$h_{ij}^{\textrm{P}}$/$\phi_{ij}$  & \tabincell{c}{Head loss on the pipe/valve from \\node $i$ to node $j$}\\
\hline
$h_{ij}^{\textrm{M}}$  & \tabincell{c}{Head gain on the pump from \\node $i$ to node $j$}\\
\hline
$q_{ij}$  & \tabincell{c}{Flow through on the link from \\node $i$ to node $j$}\\
\hline
$q_{ij}^M$/$q_{ij}^P$/$q_{ij}^V$  & \tabincell{c}{Flow through on pump/pipe/valve from \\ node $i$ to node $j$}\\
\hline
\hline
\end{tabular}
\end{center}
\end{table}

\subsection{Leveraging Linear Approximation of Hydraulic Coupling}
The nonlinear energy conservation \eqref{eq:general_energy_conservation} renders the water flow formulation nonconvex. This hinders the development of a computationally efficient stochastic optimal water flow problem where distributionally robust optimization and risk measures can be utilized to balance system performance and robustness. To that end, we provide a Laplacian-based linearization that utilizes the successive linearization algorithm to enable a highly accurate approximation of the original nonlinear energy conservation \eqref{eq:general_energy_conservation}.

The energy conservation \eqref{eq:pump_head_gain}, \eqref{eq:pipe_head_loss} and \eqref{eq:valve_head_reduce} can be concluded in a compact matrix form
\begin{equation}\label{eq:matrix_energy_conservation}
    \mathbf{B}_f \mathbf{h}(t) = \mathbf{q}^\intercal(t) \mathbf{N} \mathbf{q}(t) + \mathbf{P}\mathbf{q}(t) + \mathbf{q}_0,
\end{equation}
where $\mathbf{h}(t):=[h_1(t),\ldots, h_N(t)]^\intercal \in \mathbf{R}^{N}$ and $\mathbf{q}(t):= \{q_{ij}(t)| (i,j) \in \mathcal{I}\cup \mathcal{M}\} \cup \{\phi_{ij}(t)|(i,j)\in\mathcal{V}\} \in \mathbf{R}^{|\mathcal{E}|} $ collect network states, e.g., head, flow and valve settings. The constant matrices/vector $\mathbf{N}\in\mathbf{R}^{|\mathcal{E}| \times|\mathcal{E}|}$, $\mathbf{P}\in\mathbf{R}^{|\mathcal{E}| \times|\mathcal{E}|}$ and $\mathbf{q}_0\in\mathbf{R}^{|\mathcal{E}|}$ explicitly depend on the property of pipelines, pump and valves. The incident matrix of graph $\mathcal{G}$ is denoted by $B_f \in \mathbf{R}^{|\mathcal{E}|\times|N|}$, having entries
\begin{equation}\label{eq:incidence_matrix}
\begin{aligned}
& \mathbf{B}_f(n,i) \\
& = \left\{ \begin{array}{ll}
1 & \textrm{if flow in $n^\mathrm{th}$ link is away from $i^\mathrm{th}$ node }\\
-1 & \textrm{if flow in $n^\mathrm{th}$ link is towards $i^\mathrm{th}$ node } \\
0 &  \textrm{if flow in $n^\mathrm{th}$ link is not incident on $i^\mathrm{th}$ node }
\end{array} \right ..
\end{aligned}
\end{equation}

The nonlinearities in \eqref{eq:matrix_energy_conservation} make the OWF problem non-convex and computationally challenging. Therefore, in the rest of this subsection, we will seek to linearize \eqref{eq:matrix_energy_conservation} instead. We express the flow as $\mathbf{q} = \mathbf{\bar{q}} + \Delta \mathbf{q}$, where $\bar{\mathbf{q}}\in\mathbf{R}^{|\mathcal{E}|}$ is the nominal water flow vector and $\Delta \mathbf{q}\in\mathbf{R}^{|\mathcal{E}|}$ captures disturbances around the nominal values. To lighten notation we omit the time index in the discussion of linearization in this section. Substituting $\mathbf{q} = \mathbf{\bar{q}} + \Delta \mathbf{q}$ into \eqref{eq:matrix_energy_conservation}, we have
\begin{equation*}\label{eq:linearization}
\begin{aligned}
    \mathbf{B}_f\mathbf{h}  & = \left(\mathbf{\bar{q}} + \Delta \mathbf{q}\right)^\intercal \mathbf{N} \left(\mathbf{\bar{q}} + \Delta \mathbf{q}\right) + \mathbf{P}\left(\mathbf{\bar{q}} + \Delta \mathbf{q}\right) + \mathbf{q}_0 \\
    & + \mathbf{\bar{q}}^\intercal \mathbf{N}\mathbf{\bar{q}} + \mathbf{P}\mathbf{\bar{q}} + \left(2\mathbf{\bar{q}}^\intercal \mathbf{N} + \mathbf{P} \right)\Delta \mathbf{q} + \Delta\mathbf{q}^\intercal \mathbf{N} \Delta\mathbf{q} + \mathbf{q}_0.
\end{aligned}
\end{equation*}
Neglecting second-order terms in $\Delta \mathbf{q}$, \eqref{eq:matrix_energy_conservation} becomes approximately
\begin{equation}\nonumber
    \mathbf{B}_f\mathbf{h} \approx  \underbrace{ \bar{\mathbf{q}}^\intercal \mathbf{N} \bar{\mathbf{q}} + \mathbf{P} \bar{\mathbf{q}}+ \mathbf{q}_0}_{\mathbf{A}} + \underbrace{\left(2\bar{\mathbf{q}}^\intercal \mathbf{N} + \mathbf{P} \right)}_{\mathbf{B}} \Delta \mathbf{q},
\end{equation}
where matrices $\mathbf{A},\mathbf{B}\in \mathbf{R}^{|\mathcal{E}|\times|\mathcal{E}|}$. Now, we turn our attention to solving for the water flow perturbation vector $\Delta \mathbf{q}$. Decomposing all energy conservation on each pipeline, we can write the above linearization in the scalar form
\begin{equation*}
    h_i - h_j = a_{ij} + b_{ij}\Delta q_{ij}, \quad \forall (i,j) \in \mathcal{E},
\end{equation*}
where $a_{ij}$ and $b_{ij}$ denote the elements of $\mathbf{A}$ and $\mathbf{B}$, respectively. The water flow perturbation on each link is
\begin{equation*}
    \Delta q_{ij} = \frac{1}{b_{ij}}(h_i-h_j) - \frac{1}{b_{ij}}a_{ij}, \quad \forall (i,j)\in \mathcal{E}.
\end{equation*}
The sum of water perturbation carrying away from $i^\mathrm{th}$ node is defined as $\Delta Q_i\in \mathbf{R}$ around the nominal operation point $\mathbf{\bar{q}}$ given by
\begin{equation}\label{eq:perturbation_flow_leaving}
\Delta Q_i = \sum_j \left[\frac{1}{b_{ij}}(h_i - h_j)\right] - \underbrace{\sum_j \left[\frac{a_{ij}}{b_{ij}}\right]}_{\bar{Q}_i}, \quad \forall i \in \mathcal{N}.
\end{equation}
Note that the first term in \eqref{eq:perturbation_flow_leaving} can be expressed using the network Laplacian matrix $\mathbf{L}\in \mathbf{R}^{N \times N}$ defined by the edge weights $\frac{1}{b_{ij}}$. The second term in \eqref{eq:perturbation_flow_leaving} is the nominal carrying flow $\bar{Q}_i\in\mathbf{R}$ at $i^\mathrm{th}$ node. Then the linear energy conservation of water distribution network in compact form is  
\begin{equation*}
    \mathbf{L}\mathbf{h} = \Delta\mathbf{Q} + \mathbf{\bar{Q}},
\end{equation*}
where the network's Laplacian matrix $\mathbf{L}$ has elements 
\begin{equation}\label{eq:laplacian_matrix}\nonumber
\mathbf{L}_{ij} = \left\{ \begin{array}{ll}
\sum_{l \sim i} \frac{1}{b_{il}}  & \textrm{if $i=j$}\\
-\frac{1}{b_{ij}} & (i,j) \in \mathcal{E} \\
0 &  (i,j) \notin \mathcal{E}
\end{array} \right ..
\end{equation}
We define two vectors as $\Delta \mathbf{Q} := [\Delta Q_1, \ldots,\Delta Q_N]^\intercal$ and $\mathbf{\bar{Q}}:=[ \{\bar{Q}_i\in\mathbf{R} | \bar{Q}_i=\sum_j \frac{a_{ij}}{b_{ij}} \}$.
Given the incidence matrix defined in \eqref{eq:incidence_matrix}, the following energy conservation constraint holds as a function of water flow perturbations
\begin{equation}\label{eq:Laplacian_Energy_Conservation}
    \mathbf{L}\mathbf{h} = \mathbf{B}_f^\intercal \Delta \mathbf{q} + \mathbf{\bar{Q}},
\end{equation}
where the Laplacian-based compact form maps the water flow disturbance $\Delta\mathbf{q}$ to the network head $\mathbf{h}$ around the linearized point of nodal water carrying $\mathbf{\bar{Q}}$.

\subsection{Verifying Feasibility of Laplacian Approximation}
To validate the effectiveness and feasibility of the proposed Laplacian approximation, we solve a water flow feasibility problem
\begin{equation}\label{eq:feasible problem}
\begin{aligned}
    \textbf{WFP-0:} ~~~~&\min_{\mathbf{h},\Delta \mathbf{q}} & 0~~~~~~~~~ \\
    ~ & \textrm{subject to} & \eqref{eq:mass_balance} ~\textrm{and} ~\eqref{eq:Laplacian_Energy_Conservation}.
\end{aligned}
\end{equation}
by utilizing successive linearization algorithm and compare the solutions to the water flow results from EPANET \cite{rossman2000epanet} modelled via the nonlinear energy conservation constraints \eqref{eq:matrix_energy_conservation}. The overall successive linearization process is presented in Algorithm \ref{algorithm:successive}. We use $\delta\in\mathbf{R}_{++}$ to indicate the stop criteria of Algorithm 1, which limits the acceptable discrepancy between two successive optimization results. It can be set to a small positive constant.

We empirically observe that the successive linearization algorithm for WDN in Algorithm \ref{algorithm:successive} provides a good numerical approximation of the nonlinear hydraulic dynamics $(\mathbf{q},\mathbf{h})$ and flow directions if the following assumptions hold. More numerical results and discussions illustrating the approximation accuracy are presented in Appendix A.
\begin{assumption}
The water distribution network $\mathcal{G}$ is a tree network \cite{bullo2019lectures}.
\end{assumption}

\begin{assumption}
Pumps and valves are all active, i.e., $z_{ij} =1, \forall (i,j)\in\mathcal{M}\cup\mathcal{V}$.
\end{assumption}
The main purpose of this paper is to seek the optimal pump schedules satisfied the network constraints (e.g., flows and head limits) under network uncertainties. The proposed stochastic OWF formulation is designed for the real-time optimal control, which is based on the pre-determined operation status of all active devices, i.e., $z_{ij}$ is not an optimization variable in \eqref{eq:general_energy_conservation}. Additionally, we assume the network topology is a pure tree. This ensures that the successive linearization algorithm converges \eqref{eq:feasible problem} to a nonlinear feasible point and gives the correct directions of water flows.

\begin{algorithm}
    \caption{Successive Linearization Algorithm for Water Distribution Networks}\label{algorithm:successive}
     \hspace*{\algorithmicindent} \textbf{Input:}.EPANET Network Information \texttt{.inp} source file\footnotemark, demand $\mathbf{d}$, incidence matrix $B_f$ \\
    \hspace*{\algorithmicindent} \textbf{Output:} System operating points $\mathbf{h}^*$ and $\mathbf{q}^*$
    \begin{algorithmic}[1]
    \State Initialize: $n = 0$, initial nominal water flow $\mathbf{\bar{q}}_0$ and Linearization error Err $> \delta$
        \While{$\textrm{Err}>\delta$}
        \State Calculate Laplacian $L_n$ and flow carrying vector $\bar{Q}_n$ for $n^\mathrm{th}$ iteration
        \State Solve \textbf{WFP-0} for $\Delta \mathbf{q}_n^*$ and $\mathbf{h}_n^*$ from \eqref{eq:feasible problem}
        \State Compute the linearization error, $\textrm{Err} = \|\Delta \mathbf{q}_n^*\|_2^2$
        \State Update nominal water flow points $\mathbf{\bar{q}}_{n+1} = \mathbf{\bar{q}}_n + \Delta \mathbf{q}_n^*$
        \State $n = n +1$
        \EndWhile\\
        \Return solutions $\mathbf{q}^* = \mathbf{\bar{q}}_{n+1}$, $\mathbf{h}^* = \mathbf{h}_{n}$    
\end{algorithmic}
\end{algorithm} 

\footnotetext{The EPANET network information source file contains the topology of water networks and the properties of the hydraulic components (e.g., pipes, pump, tank and etc.).}

\begin{remark}(Hydraulic-Network Simplification). Our proposed framework focuses on optimal operation for water distribution networks instead of hydraulic design and analysis. Therefore, we assume all networks are trees, where the proposed linearization approach is highly effective. Many networks can be simplified or approximated as a tree network using various techniques, which facilitates a higher-level interpretation of the main network structure \cite{perelman2012water,anderson1995hydraulic}.
\end{remark}

\begin{remark}(Successive Linearization Initialization). The initial nominal water flow $\mathbf{q}_0$ is an input of successive linearization algorithm for WDNs. We suggest here a possible $\mathbf{q}_0$ for various components (e.g., pumps and pipes) to initialize Algorithm \ref{algorithm:successive} for improved convergence. For all pipes $\mathcal{I}$ in WDNs, the initial water flow is corresponding to the flow speed $1$ CFS \cite{rossman2000epanet}. The actual input initial water flow is adjusted based on the properties of individual pipeline (i.e., length and diameter). The initial linearized water flow point of pumps will come from the pump efficiency curve \cite{verleye2013optimising}

\begin{equation}\nonumber
    E_{ij} = e_{ij}^1 q^2_{ij} + e^2_{ij} q_{ij} + e^3_{ij}, \quad \forall (i,j)\in\mathcal{M}.
\end{equation}
The successive linearization of pumps starts from the most efficient point as
\begin{equation}\nonumber
    \frac{\partial E_{ij}}{\partial q_{ij}} = 2 e^1_{ij} q_{ij} + e^2_{ij} = 0, \quad \bar{q}_{ij} = \frac{e^2_{ij}}{2e^1_{ij}}, \quad \forall (i,j) \in \mathcal{M}.
\end{equation}
Empirically speaking, starting points satisfying the physical constraints will lead to a feasible solution.
\end{remark}
The proposed Laplacian-based approximation of hydraulic flow-head coupling provides a basic model of water distribution networks for us to develop the data-based distributionally robust optimal pump control. In the rest of the paper, we first present the general idea and then leverage the Laplacian-based successive linearization to promote a computationally-efficient framework for optimal pump control and nodal pressure management.

\section{Data-based Multi-period Stochastic Optimal Water Flow}\label{sec:stochastic_OWF}
In this section, we formulate a stochastic OWF problem as a distributionally robust stochastic optimal control problem. We first pose the problem generically to highlight the overall approach and in subsequent sections we incorporate the linearization of hydraulic modelling in Section \ref{sec:mode and linearization} for a tractable and computationally-efficient stochastic OWF. This framework is more general than most stochastic OWF in the literature, which typically focus only on individual or single-stage optimization problems, or has a less sophisticated approach for explicitly incorporating uncertainties. Consider a multi-period data-driven distributionally robust optimization problem
\begin{subequations} \label{eq:general_stochastic_OWF}
\begin{align}
	&\inf_{\pi \in \Pi} \sup_{\mathds P \in \mathcal P} && \mathds{E}^{\mathds{P}} \sum_{t = 0}^{T} h_t(x_t, u_t, \xi_t), \\
	& \subjectto && x_{t + 1}  = f_t (x_t, u_t, \xi_t),\\
	& &&u_t = ~ \pi(x_0, \ldots, x_t, \bm{\xi}_t, \mathcal{D}_t),\\
	& && (x_t, u_t) \in \mathcal{X}_t,
\end{align}
\end{subequations}
where $x_t \in \mathbf{R}^n$ represent the state vector at time $t$ that includes the internal states of all elements (i.e., valves, tanks and pipes). Let $u_t \in \mathbf{R}^m$ denote a control input vector that includes inputs for all controllable components (e.g., pump output and valve settings). The $\xi_t \in \mathbf{R}^{N_\xi}$ denote random vectors in a probability space $(\Omega, \mathcal{F}, \mathds{P}_t)$ which includes forecast errors of all uncertainties in the network. 

The goal of \eqref{eq:general_stochastic_OWF} is to find a optimal feedback policy that minimizes the expected value of the system objective function $h_t: \mathbf{R}^n \times \mathbf{R}^m \times \mathbf{R}^{N_\xi} \rightarrow \mathbf{R}$ robust to the worst-case distribution in the forecast error ambiguity set $\mathcal{P}$. We consider a setting where the objective function $h_t$ includes both operating costs and risks of violating various network and device constraints and is assumed to be continuous and convex as functions of $(x_t, u_t)$ for any fixed $\xi_t$. The system dynamics function  $f_t: \mathbf{R}^n \times \mathbf{R}^m \times \mathbf{R}^{N_\xi} \to \mathbf{R}^n$ models internal dynamics of all network-connected components, such as water storage tanks. The general feasible set $\mathcal{X}_t$ includes other network and device constraints, such as mass balance, energy conservation, operational bounds on nodal heads and pipe flows (some constraints may be modeled deterministically with respect to mean values and others may be included as risk terms in the objective function).

Since the real distributions of forecast errors are never known in practice, we explicitly account for uncertainty in their distributions themselves by assuming that the real but unknown distribution $\mathds{P}_t$ belongs to an \textit{ambiguity set} $\mathcal{P}_t$ of distributions which will be constructed from a forecast sampling dataset.We collect the forecast error over an operating horizon $t$ as $\bm{\xi}_t := [\xi_1^\intercal,\ldots,\xi_{t}^\intercal]^\intercal \in \mathbf{R}^{N_\xi t}$, which has joint distribution $\mathds{P}$ and corresponding ambiguity set $\mathcal{P}$. 

In this multi-period stochastic OWF, we are seeking a series of closed-loop feedback policies in the form $u_t = \pi(x_0, \ldots, x_t, \bm{\xi}_{0:t}, \mathcal{D}_t)$ explicitly considering forecast errors describing historical patterns, where the term $\mathcal{D}_t$ indicates all network component model information and the parameterization of the ambiguity set of the forecast error distribution. This framework allows for design of not only for current nominal reaction, but also reactions to future uncertainty realizations. The policy function $\pi$ maps all available information to control actions and is an element of a set $\Pi$ of measurable functions.

\subsection{Ambiguity Sets based on Wasserstein Metric}
One of the main challenges for solving \eqref{eq:general_stochastic_OWF} is how to utilize our available information of uncertainties to appropriately realize the distributions for a tractable problem reformulation. There is a variety of ways to reformulate the general stochastic OWF problem \eqref{eq:general_stochastic_OWF} to obtain tractable subproblems that can be solved by standard convex optimization solvers. These include assuming specific functional forms for the forecast error distribution (e.g., Gaussian) \cite{sampathirao2018gpu} and using specific constraint risk functions, such as those encoding value at risk (i.e., chance constraints) \cite{grosso2014chance,grosso2017stochastic}, conditional value at risk (CVaR) \cite{guo2020distributionally}, distributional robustness \cite{guo2020distributionally} and support robustness \cite{goryashko2014robust}. In all cases, the out-of-sample performance of the resulting decisions in operational practice ultimately relies on  1) the quality of data describing the forecast errors and 2) the validity of assumptions made about probability distributions. Many existing approaches make either too strong or too weak assumptions that possibly lead to underestimation or overestimation of risk.

In this paper, we utilize a recently proposed tractable method \cite{mohajerin} in a multi-period data-based stochastic OWF, in which the ambiguity set is based on a finite forecast error training dataset $\hat{\Xi}^{N_s}$ via Wasserstein balls. Comparing with others existing ambiguity sets \cite{li2017ambiguous,wiesemann2014distributionally,bertsimas2017robust,wozabal2012framework,jiang2016data,erdougan2006ambiguous}, Wasserstein balls offer the powerful out-of-sample performance and allow water distribution network operators to control the conservativeness of the decisions, which promote the flexibility of water distribution network from a practical perspective. We optimize an expected objective over the worst-case distribution in the ambiguity set $\mathcal{P}$, which can be formulated as a finite-dimensional convex program. The decisions from this stochastic OWF provide an upper confidence bound under forecast errors realization, quantified by the size of the \textit{ambiguity set} (i.e., Wasserstein radius \cite{mohajerin}).
The Wasserstein metric defines a distance in the space $\mathcal{M}(\Xi)$ of all probability distributions $\mathds{Q}$ supported on a set $\Xi$ with $\mathds{E}^{\mathds{Q}}[\| \xi \|] = \int_{\Xi} \| \xi \| \mathds{Q}(d\xi) < \infty$.

\begin{definition}[Wasserstein Metric \cite{kantorovich1958space},\cite{esfahani2018data}] Given all distributions $\mathds{Q}_1$,$\mathds{Q}_2$ supported on $\Xi$, the Wasserstein metric $d_W:\mathcal{M}(\Xi)\times \mathcal{M}(\Xi) \to \mathbf{R}_+$ is defined as
\begin{equation}\nonumber
    d_W\left(\mathds{Q}_1,\mathds{Q}_2\right):=\int_{\Xi}\|\xi_1 - \xi_2\| \Theta(d\xi_1,d\xi_2),
\end{equation}
where $\Theta$ represent a joint distribution of $\xi_1$ and $\xi_2$ with marginals $\mathds{Q}_1$ and $\mathds{Q}_2$, respectively and $\|\cdot\|$ indicates an arbitrary norm on $\mathbf{R}^{N_{\xi}}$.
\end{definition}

The Wasserstein metric quantifies the ``transportation costs" to move mass from one distribution to another. The Wasserstein ambiguity set is defined by
\begin{equation} \label{ambiguityset}
\hat{\mathcal{P}}^{N_s} := \bigg\{ \mathds{Q} \in \mathcal{M}(\Xi): d_W(\hat{\mathds{P}}^{N_s}, \mathds{Q}) \leq \epsilon \bigg\}.
\end{equation}
This ambiguity set $\hat{\mathcal{P}}^{N_s}$ constructs a ball with radius $\epsilon$ in Wasserstein distance around the empirical distribution $\hat{\mathds{P}}^{N_s}$ on the training dataset. The radius $\epsilon$ can be chosen so that the ball contains the true distribution $\mathds{P}$ with a prescribed confidence level and leads to performance guarantees \cite{mohajerin}. The radius $\epsilon$ also explicitly controls the conservativeness of the resulting decision. Large $\epsilon$ would produce decisions that rely less on the specific features of the uniform empirical distribution supported by the training dataset $\hat \Xi^{N_s}$ and improve robustness to inherent sampling errors. We will discuss the use of this conservativeness index for our stochastic OWF problem.

\subsection{Data-based Distributionally Robust Model Predictive Control of Optimal Water Flow}
The goal of our data-based distributionally robust stochastic OWF framework is to interpret and demonstrate inherent trade-offs between efficiency and risk of constraint violations. Accordingly, the objective function comprises a weighted sum of an operational cost function and a constraint violation risk function: $h_t =  J_\textrm{Cost}^t + \rho J_\textrm{Risk}^t$, where $\rho \in \mathbf{R}_+$ is a weight that quantifies the network operator's risk aversion. The operational cost function is assumed to be linear or convex quadratic. The cost functions will be discussed in detail in Section \ref{sec:data-driven_OWF}.

The constraint violation risk function $J_\textrm{Risk}$ comprises a sum of the conditional value-at-risk (CVaR) \cite{rockafellar2000optimization} of a set of $N_\ell$ network and device constraint functions. The conditional value-at-risk is a well known and coherent risk measurement in finance \cite{rockafellar2000optimization}. Here we introduce the CVaR risk metric to solve a MPC-based OWF engineering problem, due to the large variation of water demand uncertainties. Minimizing the CVaR of constraint violation limits both the frequency and expected severity of constraints. Specifically, we have
\begin{equation}\label{generalDRO}\nonumber
J_\textrm{Risk}^t := \sum_{i=1}^{N_\ell} \text{CVaR}_{\mathds{P}}^\beta [\ell_i(x_t, u_t, \xi_t) ],
\end{equation}
where $\beta \in (0,1]$ refers to the confidence level of CVaR under the distribution $\mathds{P}$ of random variable $\xi_t$. The risk measure encoded by the CVaR metric for the random variable $\xi$ at level $\beta$ is stated as $\textrm{CVaR}_\mathbb{P}^\beta(\ell_i(x_t,u_t,\xi_t)) := \inf_{\kappa_i^t}\mathbb{E}\frac{1}{\beta} [\ell_i(x_t,u_t,\xi_t) - \kappa_i^t]_+ + \kappa, \forall i, t$, where $\kappa_i^t \in\mathbf{R}_+$ is an auxiliary variable. Intuitively, the constraint violation risk function $J_{\textrm{Risk}}$ could be understood as the sum of networks and devices constraint violation magnitude at risk level $\beta$. The details of CVaR constraint convex reformulation are shown in the next Section.

The general problem \eqref{eq:general_stochastic_OWF} will be approached with a distributionally robust model predictive control (MPC) algorithm. MPC is a feedback control technique that solves a sequence of open-loop optimization problems over a planning horizon $\mathcal{H}_t$ (which in general may be smaller than the overall horizon $T$). At each time $t$, we solve the distributionally robust optimization problem over a set $\Pi_{\text{affine}}$ of affine feedback policies using the Wasserstein ambiguity set \eqref{ambiguityset}, where $\pi_{\textrm{affine}}$ collects linear functions which map the past system states, uncertainties and historical data to control actions.\\
\textbf{Distributionally Robust MPC for Stochastic OWF:} 
\begin{subequations} 
\label{general_DROOWF_rho}
\begin{align}
	&\inf_{\pi \in \Pi_{\text{affine}}} \sup_{\mathds P \in \hat{\mathcal{P}}^{N_s} } && \mathds{E}^{\mathds{P}} \sum_{\tau = t}^{t+\mathcal{H}_t} J^\tau_{\text{Cost}} + \rho J^\tau_{\text{Risk}}, \\
	& \subjectto && x_{\tau + 1}  = f_\tau (x_\tau, u_\tau, \xi_\tau), \\
	& &&u_\tau = ~ \pi(x_0, \ldots, x_\tau, \bm{\xi}_\tau, \mathcal{D}_\tau), \\
	& && (x_\tau, u_\tau) \in \mathcal{X}_\tau.
\end{align}
\end{subequations}
Only the immediate control decisions for time $t$ are implemented on the controllable device inputs. Then time shifts forward one step, new forecast errors and states are realized, the optimization problem \eqref{general_DROOWF_rho} is re-solved at time $t+1$ and the process repeats. This approach allows any forecasting methodology to be utilized to predict uncertainties over the planning horizon. Furthermore, the forecast error dataset $\hat{\mathds{P}}^{N_s}$, which defines the center of the ambiguity set $\hat{\mathcal{P}}^{N_s}$, can be updated online as more forecast error data is obtained. It is also possible to remove outdated data online to account for time-varying distributions.

In the rest of the paper, we will use the specific model of water distribution networks discussed in Section \ref{sec:mode and linearization}, where the subproblems \eqref{general_DROOWF_rho} have exact tractable convex reformulations as quadratic programs \cite{mohajerin} and can be solved to global optimality with standard solvers.

\section{Chance-constraints and Distributionally Robustness Formulation}
\label{sec:data-driven_OWF}
Following our proposed formulation above, we begin this section by introducing the state space expression of WDN hydraulic dynamics, briefly discuss chance constraints and describe a convex reformulations of the stochastic optimal water flow problem based on conditional value-at-risk and distributionally robust optimization.

\subsection{Network Dynamics in State-Space Format}
The WDN model discussed in the previous section can be summarized in a difference algebraic equation (DAE) model
\begin{subequations}\label{eq:State_space_expression}
\begin{align}
    x(t+1) & = \bar{A} x(t) + \bar{B}_uu(t) + \bar{B}_v v(t),\label{eq:tank_SS}\\
    d(t) & = \bar{E}_u u(t) + \bar{E}_v v(t), \label{eq:mass_balance_SS}\\
    \bar{F}_x x(t) + \bar{F}_l l(t) & = \bar{F}_u u(t) + \bar{F}_v v^{\mathrm{P}}(t) + \bar{F}_\phi \phi(t) + \bar{F}_0, \label{eq:energy_conservion_SS}
\end{align}
\end{subequations}
where the decision variables $\{x,u,l,v,v^{\mathrm{P}},\phi\}$ are defined in Table 2 and the constant matrices $\{A,B,E,F\}$ are derived from the hydraulic dynamics in Section \ref{sec:mode and linearization}. We detail these constants in term of the Laplacian-based hydraulic model \eqref{eq:Laplacian_Energy_Conservation} in the Appendix B. The dynamics of tank head \eqref{eq:tank_head_dynamics} is given in \eqref{eq:tank_SS} and the mass balance \eqref{eq:mass_balance} and linearized energy conservation \eqref{eq:Laplacian_Energy_Conservation} are summarized in \eqref{eq:mass_balance_SS} and \eqref{eq:energy_conservion_SS}, respectively. For compact notation, we concatenate the states, inputs and demands over the planning horizon as
$\mathbf{x}_t = [x(1)^\intercal,\ldots,x(t)^\intercal]^\intercal \in \mathbf{R}^{n_{\mathrm{TK}}t}$, $\mathbf{u}_t = [u(0)^\intercal,\ldots,u(t-1)^\intercal]^\intercal \in \mathbf{R}^{n_ut}$, $\mathbf{v}_t = [v(1)^\intercal,\ldots,v(t)^\intercal]^\intercal \in \mathbf{R}^{n_vt}$, $\mathbf{v}^\mathrm{P}_t = [v^\mathrm{P}(1)^\intercal,\ldots,v^\mathrm{P}(t)^\intercal]^\intercal \in \mathbf{R}^{n_pt}$, $\mathbf{l}_t = [l(1)^\intercal,\ldots,l(t)^\intercal]^\intercal \in \mathbf{R}^{n_lt}$, $\bm{\phi}_t = [\phi(1)^\intercal,\ldots,\phi(t)^\intercal]^\intercal \in \mathbf{R}^{n_\phi t}$ and $\mathbf{d}_t = [d(0)^\intercal,\ldots,d(t-1)^\intercal]^\intercal\in \mathbf{R}^{Nt}$, yielding
\begin{subequations}\nonumber
\begin{align}
\mathbf{x}_t & = Ax_0 + B_u\mathbf{u}_t + B_v\mathbf{v}_t,\\
\mathbf{d}_t & = E_u\mathbf{u}_t + E_v\mathbf{v}_t,\\
F_x\mathbf{x}_t + F_l\mathbf{l}_t & = F_u\mathbf{u}_t + F_v\mathbf{v}_t^\mathrm{P} + F_\phi \bm{\phi}_t + F_0,
\end{align}
\end{subequations}
where $I_t$ indicates a $t$-dimensional identity matrix
\begin{equation}\nonumber
\begin{aligned}
& && E_u = I_t \otimes \bar{E}_u, ~~E_v = I_t \otimes \bar{E}_v, ~~F_x = I_t \otimes \bar{F}_x,\\
& && F_l = I_t \otimes \bar{F}_l, ~~F_u = I_t \otimes \bar{E}_u, ~~F_v = I_t \otimes \bar{F}_v,\\
& && F_\phi = I_t \otimes \bar{F}_\phi, ~~F_0 = I_t \otimes \bar{F}_0,\\
& &&A = \begin{bmatrix}
\bar{A}\\
\bar{A}^2\\
\vdots\\
\bar{A}^{t}
\end{bmatrix},B_d = \begin{bmatrix}
\bar{B}_u & 0 & \cdots & 0\\
\bar{A}\bar{B}_u & \bar{B}_u & \ddots & 0\\
\vdots & \ddots & \ddots & \vdots\\
\bar{A}^{t-1}\bar{B}_u & \cdots & \bar{A}\bar{B}_u & \bar{B}_u
\end{bmatrix}, \\
& &&B_v = \begin{bmatrix}
\bar{B} & 0 & \cdots & 0\\
\bar{A}\bar{B}_v & \bar{B}_v & \ddots & 0\\
\vdots & \ddots & \ddots & \vdots\\
\bar{A}^{t-1}\bar{B}_v & \cdots & \bar{A}\bar{B}_v & \bar{B}_v
\end{bmatrix}.~~~~~~~~~~~~~~~~~~~~~~~~~~~~~~~~~~~~~~~~~~~~~~~~~~~~~~~~~~~~~
\end{aligned}
\end{equation}

\begin{table}[htbp!]
\begin{center}
\caption{Variable Description in DAE Model}
\begin{tabular}{c|c|c}
\hline
\hline
N* & Description & Dimension\\
\hline
$x$  & a vector collecting heads at tanks &  $n_{\textrm{TK}} = |\mathcal{T}|$\\
\hline
$l$  & \tabincell{c}{a vector collecting heads
\\at junctions \& reservoirs} & $n_l = |\mathcal{J}|+|\mathcal{S}|$\\
\hline
$u$  & \tabincell{c}{a vector collecting flow at pumps} & $n_u = |\mathcal{M}|$\\
\hline 
$v$  & \tabincell{c}{a vector collecting flows \\through pipes \& valves} & $n_v = |\mathcal{I}|+|\mathcal{V}|$\\
\hline 
$v^{\mathrm{P}}$  & \tabincell{c}{a vector collecting flows through pipe} & $n_p = |\mathcal{I}|$\\
\hline 
$\phi$  & a vector collecting head loss on PRVs & $n_{\phi}= |\mathcal{V}|$\\ 
\hline
\hline
\end{tabular}
\end{center}
\quad N* indicates a abbreviation of Notation.
\label{table:DAE_model_definition}
\end{table}

\subsection{Cost Functions and Constraints}
Multiple objective functions can be included in the stochastic optimal water flow problem
\begin{subequations}
\begin{align}
    J_1^t & = u(t)^\intercal H_u(t) u(t) + f^\intercal_u(t) u(t) + f_0,  \label{eq:pump_objective} \\
    J_2^t & = \Delta u(t)^\intercal \Delta u(t), \label{eq: control_smooth} \\
    J_3^t & = \left(x(t) - V^{\textrm{safe}}\right)^\intercal\left(x(t) - V^{\textrm{safe}}\right), \label{eq: tank_safty}
\end{align}
\end{subequations}
where \eqref{eq:pump_objective} captures the pump operational cost based on time-varying electricity tariffs. The matrix $H_u$ is positive semi-definite. The control input variation between consecutive time slots (e.g., $\Delta u(t):= u(t) - u(t-1)$) can be also penalized in \eqref{eq: control_smooth} to avoid large transient in pipes and to satisfy treatment requirements. Additionally, tank management requires a safety head level $V^{\textrm{safe}}$ to account for unexpected demand given in \eqref{eq: tank_safty}.

The system constraints are introduced due to the physical nature of the decision variables (i.e, $\mathbf{x}$ and $\mathbf{u}$). We seek to enforce state and input constraints
\begin{subequations}
\begin{align}
\mathbf{u}^\textrm{min} \leq \mathbf{u}_t \leq \mathbf{u}^\textrm{max}, \label{eq:input_constraints}\\
\mathbf{x}^\textrm{min} \leq \mathbf{x}_t \leq \mathbf{x}^\textrm{max}, \label{eq:state_constraints}
\end{align}
\end{subequations}
where \eqref{eq:input_constraints} corresponds to actuator limits (e.g., pumps and valves) and \eqref{eq:state_constraints} captures bounds on pipe flows, nodal heads and tank levels. Here, $\mathbf{x}_\textrm{min}$ and $\mathbf{x}_\textrm{max}$ denote the minimum and maximum admissible bounds of states. The lower and upper physical limits of actuators are $\mathbf{u}_\textrm{min}$ and $\mathbf{u}_\textrm{max}$, respectively. In general, these constraints can not be violated strictly due to the mass conservation principles and physical restriction of components. For the rest of this paper, we assume that these hard bounds can be ``softened" to non-physical upper and lower bounds from a pre-specified safe operation zone, which can be violated probabilistically but results in safety or operational risk \cite{grosso2017stochastic}.
\vspace{-3mm}
\subsection{Multi-Period Stochastic Optimal Water Flow}
In a deterministic optimal water flow control problem, water demand uncertainty is not explicitly considered. Since actual water demands can exhibit large variations and unpredictability \cite{brentan2018water}, we model demand stochastically as $\mathbf{d}_t = \bar{\mathbf{d}}_t + \bm{\xi}_t$, with a nominal predicted value $\bar{\mathbf{d}} \in \mathbf{R}^{Nt}$ and a zero-mean forecast error $\bm{\xi}_t = [\xi_1^\intercal,\ldots, \xi_{t}^\intercal]^\intercal \in \mathbf{R}^{Nt}$ from a probability space $(\Omega, \mathcal{F}, \mathds{P}_\xi)$. The distribution captures spatiotemporal variations and dependencies among the demands.

To explicitly account for this stochasticity of water demands, we formulate the following general stochastic optimal water flow problem to find an optimal strategy for responding to forecast errors via an optimal control policy for the flow actuators $\mathbf{u}_t = \pi_t(\bm{\xi}_t)$, where $\pi_t : \mathbf{R}^{Nt} \to \mathbf{R}^{n_ut}$ is a function from a set $\Pi_c$ of causal policies. Specifically, we consider a multi-period optimal water flow problem with finite time horizon $T$
\begin{subequations}\label{generalstochasticOWF}
\begin{align}
    & \inf_{\pi_T\in\Pi_c}&& \sum_{\tau = 1}^{T} \mathds{E}^{\mathds{P}_\xi} \big[ J^\tau(\mathbf{x}_\tau, \pi_\tau(\bm{\xi}_\tau), \bm{\xi}_\tau \big],\label{objectiveFunction}\\
    & \subjectto && \mathbf{d}_T = E_u\pi_T(\bm{\xi}_T) + E_v\mathbf{v}_T,\label{massjunctionMultistage}\\
    & && \mathbf{x}_T = Ax_0 + B_u\pi_T(\bm{\xi}_T)+ B_v\mathbf{v}_T,\label{statedynamicsMultistage}\\
    & && F_x\mathbf{x}_T + F_l\mathbf{l}_T = F_u\mathbf{u}_T + F_v\mathbf{v}_T^\mathrm{P} + F_\phi \bm{\phi}_T + F_0,\label{energyconservationMultistage}\\
    \color{blue}
    & && \mathds{E}~ \mathcal{R} \big(\mathbf{u}_\textrm{min} 
    - \pi_T(\bm{\xi}_T) \big)\leq 0 ,\label{generalcontrolinequalityLower}\\
    & && \mathds{E}~\mathcal{R} \big(\pi_T(\bm{\xi}_T) - \mathbf{u}_\textrm{max} \big)\leq 0 , \label{generalcontrolinequalityUpper}\\
    & && \mathds{E}~\mathcal{R} \big( \mathbf{x}_\textrm{min} - \mathbf{x}_T \big)\leq 0 ,\label{generalstateinequalityLower}\\
    & && \mathbb{E}~\mathcal{R} \big( \mathbf{x}_T - \mathbf{x}_\textrm{max} \big)\leq 0 , \label{generalstateinequalityUpper} \color{black}
\end{align}
\end{subequations}
where $\mathcal{R}$  indicates a generic transformation of the inequality constraints into the stochastic versions with different uncertainty assumptions and stochastic optimization techniques. These include assuming specific probability distributions (i.e., Gaussian), using specific constraint risk measurement, such as value-at-risk and conditional value-at-risk, sample average approximation, scenario-approach, distributionally robust optimization and robust support.
Note that this transformation can be different in general for each constraint. For constraints that represent physical limits, we consider tightened non-physical upper and lower bounds on states and inputs from a pre-specified safe operation zone, which can be violated probabilistically but results in safety or operational risks \cite{grosso2017stochastic}. 
Since optimizing over general policies makes
problem \eqref{generalstochasticOWF} infinite dimensional, we optimize instead over a set of affine control policies
\begin{equation}
    \mathbf{u}_\tau = D_\tau\bm{\xi}_\tau + e_\tau,\\
\end{equation}
where $e_\tau \in \mathbf{R}^{n_u\tau}$ represents a nominal plan for pumps and the block lower-triangular matrix $D_\tau \in \mathbf{R}^{n_u\tau \times N\tau }$ ensures that the controller is causal. In this case, the input design variables turn to an uncertainty feedback matrix $D_\tau$ and nominal input vector $e_\tau$. Unlike traditional state-driven feedback control, the optimal feedback matrix $D_\tau$ acts as reserve policies of pumps to respond to realized water demand variations $\bm{\xi}_\tau$.

Substituting the affine control policies into \eqref{generalstochasticOWF}, the objective function \eqref{objectiveFunction} becomes convex quadratic in $D_\tau$ and $e_\tau$ and depends on the distributional information of $\bm{\xi}_\tau$. Since the policy is affine, the robust equality constraint \eqref{massjunctionMultistage} is equivalent to
\begin{equation}\label{affineEquality}
    E_uD_T = \mathbf{1}, ~~  \mathbf{\bar{d}}_T = E_ue_T +  E_v\mathbf{v}_T.
\end{equation}
With affine policies, \eqref{generalcontrolinequalityLower}-\eqref{generalstateinequalityUpper} become
\begin{subequations}
\begin{align}
& \mathcal{R}\left(D_T\bm{\xi}_T + e_T - \mathbf{u}_{\textrm{max}} \right) \leq 0,\label{eq:StochasticLinear_umax}\\
& \mathcal{R}\left(\mathbf{u}_{\textrm{min}} - D_T\bm{\xi}_T - e_T \right)\leq 0,\label{eq:StochasticLinear_umin}\\
& \mathcal{R}\left(Ax_0 + B_u(D_T\bm{\xi}_T + e_T) + B_v\mathbf{v}_T - \mathbf{x}_{\textrm{max}} \right) \leq 0,\label{eq:StochasticLinear_xmax}\\
& \mathcal{R}\left(\mathbf{x}_{\textrm{min}} - Ax_0 - B_u(D_T\bm{\xi}_T + e_T) - B_v\mathbf{v}_T \right)\leq 0,\label{eq:StochasticLinear_xmin}
\end{align}
\end{subequations}
We collect all above affine constraints inside the risk measures \eqref{eq:StochasticLinear_umax}--\eqref{eq:StochasticLinear_xmin} into a set $\mathbb{V}_{\{1:T\}}$ of $N_\ell = 2T(n_{\mathrm{TK}} + n_u)$ constraints and the expressions inside the brackets can be written in a general linear form $a_i(D_T)^\intercal\bm{\xi}_T + b_i(e_T)$, where index $i$ refers to each individual constraint in $\mathbb{V}_{\{1:T\}}$.
\subsection{Chance-Constraints}
Using a Value-at-Risk measure, the OWF problem can be posed as a chance-constrained optimization problem
\begin{subequations}
    \begin{align}
    & \inf_{D,e} && \sum_{\tau = 1}^T \mathds{E}^{\mathds{P}_\xi} \left[ J^\tau(\mathbf{x}_\tau,\mathbf{u}_\tau,\bm{\xi}_\tau)\right],\nonumber \\
    & \subjectto && \mathds{P}_\xi\left(a_i(D_T)^\intercal\bm{\xi}_T + b_i(e_T) \leq 0 \right) \geq 1 - \beta, \nonumber \\
    & && E_uD_T = \mathbf{1}, \mathbf{\bar{d}}_T = E_ue_T +  E_v\mathbf{v}_T, \nonumber  \\
    & && F_x\mathbf{x}_T + F_l\mathbf{l}_T = F_ue_T + F_v\mathbf{v}_T^\mathrm{P} + F_\phi \bm{\phi}_T + F_0, \nonumber  \\
    & && \forall i \in \mathbb{V}_{\{1:T\}}, \nonumber 
    \end{align}
\end{subequations}
where $\beta \in \mathbf{R}$ is the prescribed safety parameter or ``risk budget" for the linear constraint in set $\mathbb{V}_{\{1:T\}}$.
The subscript $\{1:T\}$ of set $\mathbb{V}_{\{1:T\}}$ indicates the set exclusively includes the state and input constraints between time interval $[1, T]$. If $\bm{\xi}_T$ is normally distributed, then it is known that the chance constraint can be written as a second-order cone constraint \cite{delage2010distributionally,calafiore2006distributionally}. However, in general chance constraints only restrict the frequency of constraint violations, not the severity. Since the real distribution is never known in practice, this approach can lead to underestimation of actual risks and poor out-of-sample performance. In this paper, we leverage a data-driven distributionally robust optimization methodology to account for both frequency and severity of constraint violation via conditional value-at-risk (CVaR) metric without assuming a particular distribution.

\subsection{Stochastic OWF based on Distributionally Robust Optimization and Conditional Value-at-Risk (CVaR)}
We treat the constraints \eqref{eq:StochasticLinear_umax}--\eqref{eq:StochasticLinear_xmin} with a risk measure derived from distributionally robust optimization techniques. It is possible to allow some constraints to be reformulated by other risk measures and optimization techniques, such as sample average approximation, moment-based distributionally robust optimization, robust optimization and Gaussian-based chance constraints. We restrict the model here only for Wasserstein metric distributionally robust techniques and leave potential combinations for the future work.

For simplicity, we consider the risk of each constraint individually; it is possible to consider risk of joint constraint violations, but this is more difficult and we leave it for future work. Recall each individual affine constraint between the finite time horizon $\mathcal{H}_t$ in the set $\mathbb{V}_{\{\underline{t}:\bar{t}\}}$ can be written in a compact form as follows. The $[\underline{t},\bar{t}]$ here refers to the finite time horizon $[t,t+\mathcal{H}_t]$.
\begin{equation*}
\begin{split}
\mathcal{C}_i^t(D_{t},e_t,\bm{\xi}_{t}) = a_i(D_{t})^\intercal\bm{\xi}_{t} + b_i(e_{t}),  t \in  [\underline{t},\bar{t}],
\end{split}
\end{equation*}
where $\mathcal{C}_i^t(\cdot)$ is the $i^\mathrm{th}$ affine constraint in the set $\mathbb{V}_{\{\underline{t},\bar{t}\}}$. The CVaR with risk level $\beta$ of the each individual constraint in the set $\mathbb{V}_{\{\underline{t},\bar{t}\}}$ is
\begin{equation} \label{defininationofCVaRdistribution}
\inf_{\kappa_i^t} \mathds{E}_{\xi_t}\bigg\{[\mathcal{C}^t_i(D_{t},e_{t},\bm{\xi}_{t}) + \kappa_i^t]_+ - \kappa_i^t\beta \bigg\} \leq 0, t \in [\underline{t},\bar{t}],
\end{equation}
where $\kappa_i^t \in \mathbf{R}$ is an auxiliary variable \cite{rockafellar2000optimization}. The expression inside the expectation in \eqref{defininationofCVaRdistribution} can be expressed in the form with risk level $\beta$
\begin{equation}\label{distributionCVaRConvex}\nonumber
\mathcal{Q}_i^t = \max_{k=1,2} \bigg[ \langle \mathbf{a}_{ik}^\beta(D_t,e_t), \xi_t \rangle + \mathbf{b}_{ik}^\beta(\kappa_i^t) \bigg], t \in [\underline{t}, \bar{t}].
\end{equation}
This expression is convex in $(D_t, e_t)$ for each fixed $\bm{\xi}_t$ since it is the maximum of two affine functions. Our risk objective function is expressed by the distributionally robust optimization of CVaR
\begin{equation}\nonumber
\begin{split}
& \hat{J}_\textrm{Risk}^t = \\
& \sum_{\tau = t}^{t+\mathcal{H}_t}\sum_{i=1}^{N_\ell} \sup_{\mathds{Q}_\tau \in \hat{\mathcal{P}}_\tau^{N_s}}\mathds{E}^{\mathds{Q}_\tau} \max_{k=1,2} \bigg[ \langle \mathbf{a}_{ik}^\beta (D_\tau, e_\tau), \hat{\bm{\xi}}_\tau \rangle + \mathbf{b}_{ik}^\beta(\kappa_i^\tau) \bigg].
\end{split}
\end{equation}
The above multi-period distributionally robust optimization can be equivalently reformulated the following quadratic program, the details of which are described in \cite{mohajerin}.
The objective is to minimize a weighted sum of an operational cost function and the total worst-case CVaR of the affine constraints in set $\mathbb{V}_{\{\underline{t},\bar{t}\}}$ (e.g., nodal head and tank level).\\
\textbf{Data-based Distributionally Robust MPC Stochastic OWF:}
\begin{subequations} \label{eq:DRO_OWF_MPC}
\begin{align}
	&\hspace{-2mm}\inf_{\begin{subarray}{c}D_\tau, e_\tau \kappa_i^\tau\\ \mathbf{v},\mathbf{v}^\mathrm{P}, \mathbf{x}, \mathbf{l}, \bm{\phi}
	\end{subarray} } \sum_{\tau = t}^{t + \mathcal{H}_t}  \bigg\{\mathds{E} [\hat{J}^\tau_{\text{Cost}}] +  \rho \hspace{-2 mm} \sup_{\mathds{Q}_\tau \in \hat{\mathcal{P}}^{N_s}_\tau} \sum_{i=1}^{N_\ell} \mathds{E}^{\mathds{Q}_\tau}[\mathcal{Q}_i^\tau]  \bigg\},  \nonumber\\
	&\hspace{-2mm} = \inf_{\begin{subarray}{c}D_\tau, e_\tau \kappa_i^\tau, \\  \lambda_i^\tau, s_{io}^\tau, \varsigma_{iko}^\tau\\\mathbf{v},\mathbf{v}^\mathrm{P}, \mathbf{x}, \mathbf{l}, \bm{\phi} \end{subarray}} \sum_{\tau=t}^{t + \mathcal{H}_t} \bigg\{\mathds{E}[\hat{J}^t_{\text{Cost}}] + \sum_{i=1}^{N_\ell} \bigg(\lambda_i\epsilon_\tau + \frac{1}{N_s}\sum_{o=1}^{N_s} s_{io}^\tau\bigg) \bigg\},\\
	&\subjectto~~~~~~~~~~~~~~~~~~~~~~~~~~~~~~~~~~~~~~~~~~~~~~~~~~~~~~~~~~~~~~~~~~~~~~~~~~\nonumber \\
	&\Big[E_uD_{\bar{t}} - \mathbf{1} \Big]_{[\underline{t},\bar{t}]}= \mathbf{0}_{N\mathcal{H}_t}, \\
	&\Big[E_ue_{\bar{t}} +  E_v\mathbf{v}_{\bar{t}} - \mathbf{\bar{d}}_{\bar{t}} \Big]_{[\underline{t},\bar{t}]} = \mathbf{0}_{N\mathcal{H}_t}, \\
	&\Big[F_x\mathbf{x}_{\bar{t}} + F_l\mathbf{l}_{\bar{t}} - F_ue_{\bar{t}} - F_v\mathbf{v}_{\bar{t}}^\mathrm{P} - F_\phi \bm{\phi}_{\bar{t}} - F_0\Big]_{[\underline{t},\bar{t}]} = \mathbf{0}_{N\mathcal{H}_t},\\
	&\rho (\langle \mathbf{a}_{ik}^\beta(D_\tau,e_\tau), \hat{\bm{\xi}}_\tau^{o}\rangle + \mathbf{b}_{ik}^\beta(\kappa_i^\tau) + \langle\varsigma_{iko}, \mathbf{z}_\tau-\mathbf{F}_\tau\hat{\bm{\xi}}_\tau^o\rangle) \le s_{io}^\tau, \label{eq:DRO_Constraint_1} \\
	&\|\mathbf{F}^\intercal_\tau \varsigma_{iko}-\rho \mathbf{a}_{ik}^\beta(D_\tau,e_\tau)\|_\infty \le \lambda_i^\tau, \label{eq:DRO_Constraint_2} \\
	&\varsigma_{iko} \geq 0, \label{eq:DRO_Constraint_3}\\
	&\forall o\le N_s, \forall i \le N_\ell, k=1,2, \tau = t,..., t + \mathcal{H}_t, \nonumber 
\end{align}
\end{subequations}
where $\rho \in \mathbf{R}_+$ quantifies the water network operators' risk aversion. This is a quadratic program that explicitly uses the training dataset $\hat{\Xi}^{N_s}_\tau = \{\hat{\bm{\xi}}_\tau^o\}_{o \leq N_s}$. The risk aversion parameter $\rho$ and the Wasserstein radius $\epsilon_\tau$ allow us to explicitly balance trade-offs between efficiency, risk and sampling errors inherent in $\hat{\Xi}^{N_s}_\tau$. The support is modeled as a polytope $\Xi_\tau := \{\bm{\xi}_\tau \in \mathbf{R}^{N_\xi\tau}: \mathbf{F}_\tau\bm{\xi}_\tau \leq \mathbf{z}_\tau \}$. The constraint $\varsigma_{iko} > 0$ holds since the uncertainty set is not-empty; on the other hand, in a case with no uncertainty (i.e, $\varsigma_{iko} = 0$), the variable $\lambda$  does not play any role and $s_{io}^\tau = \rho (\langle \mathbf{a}_{ik}^\beta(D_\tau, e_
\tau), \bm{\hat{\xi}}_\tau^{o}\rangle + \mathbf{b}_{ik}^\beta(\kappa_i^\tau))$.
\begin{figure}[!t]
    \centering
    \includegraphics[width=3.2in]{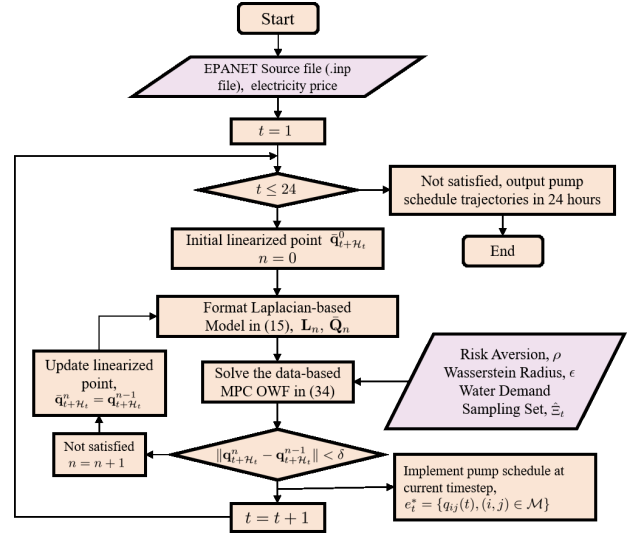}
    \caption{Flowchart of data-based distributionally robust stochastic OWF.}
    \label{fig:flowchart_DROOWF}
\end{figure}
\begin{remark}
There are three important tuning parameters in our proposed multi-period data-based stochastic OWF \eqref{eq:DRO_OWF_MPC} corresponding to different performance-risk trade-offs, which all function in the ways with their unique interpretation 
\begin{itemize}
    \item The \textbf{Wasserstein radius} $\epsilon$ improves out-of-sample performance and mitigates the effects of inherent sampling errors, which here is our main focus. The decisions optimize performance under the worst-case distributions within Wasserstein distance $\epsilon$ of the empirical distribution in probability distribution space. A larger $\epsilon$ indicates less reliance on the specific training dataset $\hat{\Xi}$ that describes the real unknown data-generating distribution, which results in more conservative decisions. The superior out-of-sample performance is achieved by this adjustable Wasserstein metric, as demonstrated in Section \ref{sec:case_studies}.
    \item \textbf{Risk aversion} $\rho$ trades off the operational risk and the nominal operational efficiency. The proposed stochastic OWF offers the system operators alternative strategies to run the water distribution networks under different risk levels. The decisions under various $\rho$ achieve various risk levels. Meanwhile, the out-of-sample performance under fixed risk aversion is controlled by the adjustable Wasserstein radius.
    \item \textbf{CVaR risk level} $\beta$ indicates the risk level of constraints \eqref{eq:StochasticLinear_umax}--\eqref{eq:StochasticLinear_xmin}, which trades off constraint violation magnitudes with nominal operational efficiency. 
\end{itemize}
\end{remark}
\noindent We fix the risk level parameter $\beta$ of CVaR to highlight the effects of $\epsilon$ and $\rho$ in the next section. It would also be interesting to explore the additional effects of changing $\beta$; however, we leave this discussion for future work.
\begin{remark}(Successive Linearization for Stochastic OWF \eqref{eq:DRO_OWF_MPC}). The data-based distributionally robust stochastic optimal water flow \eqref{eq:DRO_OWF_MPC} at $t^\mathrm{th}$ time interval is solved via successive linearization algorithm discussed in Section \ref{sec:mode and linearization}. Since all coefficients and affine constraints $\{A, B, E, F, \mathbf{a}_{ik}^\beta.\mathbf{b}_{ik}^\beta\}$ are derived from 
Laplacian-based network format \eqref{eq:Laplacian_Energy_Conservation}, at each successive linearization iteration for certain time interval, we repeatedly obtain $\{A, B, E, F, \mathbf{a}_{ik}^\beta.\mathbf{b}_{ik}^\beta\}$ in problem \eqref{eq:DRO_OWF_MPC} based on the linearized updated $\{\mathbf{L}, \mathbf{\bar{Q}}\}$ until the linearized errors converge and move to the next time period. A flowchart of proposed data-based distributionally robust stochastic OWF is demonstrated in Fig.~\ref{fig:flowchart_DROOWF}.
\end{remark}

\section{Case Studies}\label{sec:case_studies}
We now demonstrate the effectiveness of the proposed framework with numerical experiments. We use a network model derived from a portion of the Barcelona drinking water network \cite{grosso2014chance}. There are 2 reservoirs, 4 water demands, 3 tanks, 2 pumps, 4 valves and 20 junctions, the physical properties of nodes and links are given in Tables \ref{table:nodal_info} and \ref{table:link_info}, respectively. 
The nominal water demand pattern over 24 hours shown in Fig.~\ref{fig:demand_pattern} is derived from EPANET (a standard software package for analysis of drinking water distribution systems) \cite{rossman2000epanet}. Four demands are located at nodes 8, 15, 16 and 17. Realization of demand forecasting errors are generated by evaluating the so-called persistence forecast on the EPANET demand data, which predicts the water demand at the next time step to be equal to that at the previous time step. The time-of-use (TOU) electricity price is given in Fig. \ref{fig:TOU}. Generally speaking, the proposed framework can handle inputs of the electricity prices and water demand patterns.

We placed three tanks at Node 23, Node 24 and Node 25 to accommodate the water demand uncertainties associated with the downstream nodes. The lower and upper tank level in feet are restricted to $h^{\textrm{min}}_i = 525$ and $h^{\textrm{max}}_i = 530, \forall i \in \mathcal{T}$. Due to the inherent variability of water demands, tank level constraint violations may occur. Given the forecasting error data of water demand, the numerical tests are focused on reducing potential constraint violation via proposed distributionally robust framework \eqref{eq:DRO_OWF_MPC} and minimizing the operational cost under certain risk aversion as well. To have a clear and straightforward presentation, only the lower level constraints of three tanks are modelled in distributionally robust fashion \eqref{eq:DRO_Constraint_1}--\eqref{eq:DRO_Constraint_3}. Other constraints are handled via sample average approximation (SAA) \cite{linderoth2006empirical,bertsimas2018robust} or deterministic approach, though it is easy reformulate other constraints with distributional robustness.

No bound is enforced on water demand forecast errors, which implies the parameters (i.e., $\mathbf{z}$ and $\mathbf{F}$) of polytope supported set in \eqref{eq:DRO_Constraint_1}--\eqref{eq:DRO_Constraint_2} are set to zero. The variation of forecast errors increases with the prediction horizon. The number of forecast error samples in the training data set $\hat{\Xi}_t^N$ is $N_s = 100$. The stop criteria of Algorithm 1 is set to $\delta = 10^{-6}$. The simulation takes 60 seconds or less to solve DRO OWF with finite horizon $\mathcal{H}_t = 4$ (hours) using MOSEK Solver \cite{mosek} via the MATLAB interface with CVX \cite{cvx} on a laptop with 16GM of memory and a 2.8GHz Intel Core i7. 

\begin{figure}[!t]
    \centering
    \includegraphics[width=3.5in]{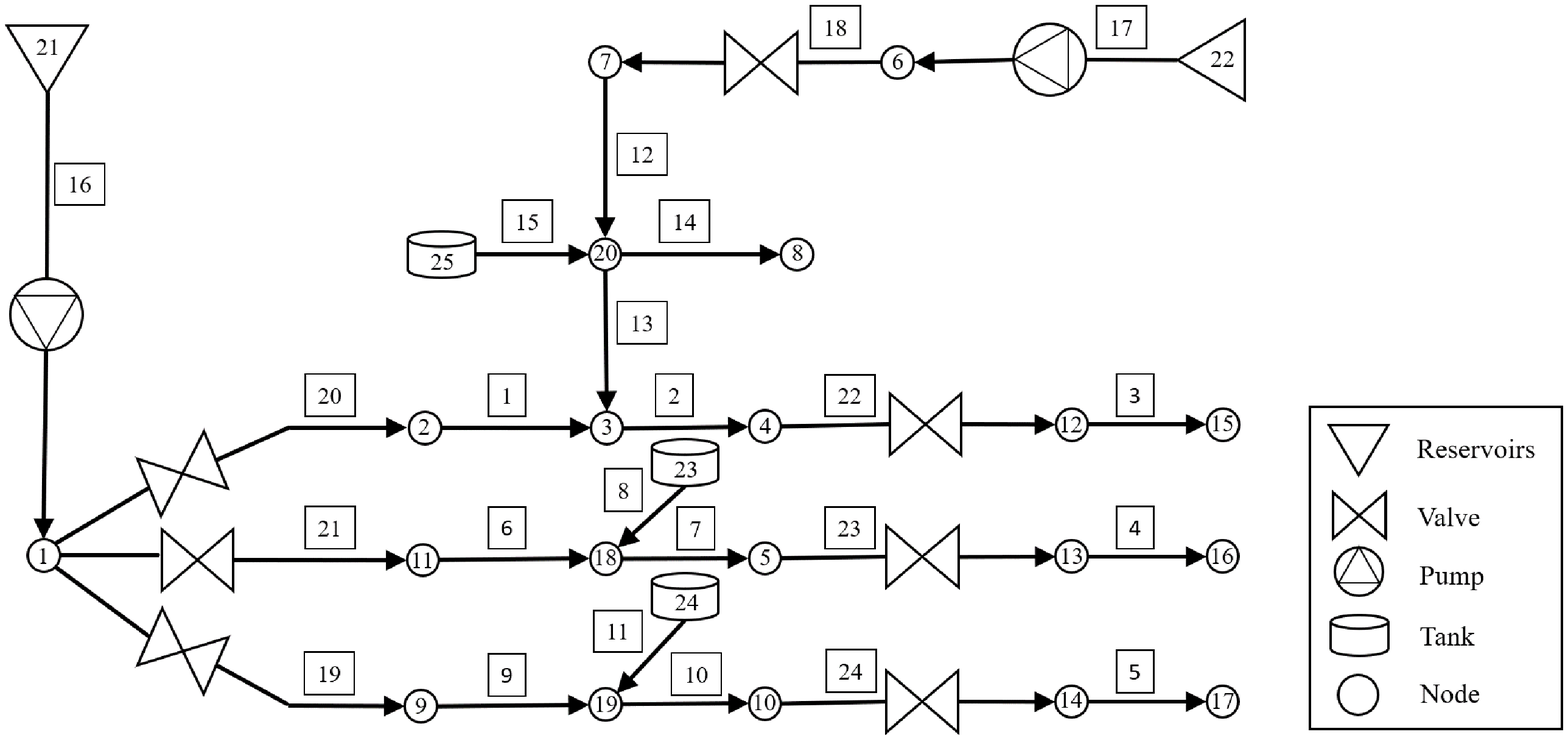}
    \caption{Barcelona drinking water network includes 25-node, 3 tanks and 2 reservoirs.}
    \label{fig:network_diagram}
\end{figure}

\begin{table}[htbp!]
\caption{Node Setting of the Barcelona Water Distribution Network}
\centering
\begin{tabular}{c|c|c|c|c|c}
\hline
\hline
Node & Type & \tabincell{c}{Base \\Demand\\(GPM)} & Node & Type & \tabincell{c}{Base \\Demand\\(GPM)}\\
\hline
1 & Junction & 0 & 2 & Junction & 0\\
\hline
3 & Junction & 0 & 4 & Junction & 0\\
\hline
5 & Junction & 0 & 6 & Junction & 0\\
\hline
7 & Junction & 0 & 8 & Junction & 100\\
\hline
9 & Junction & 0 & 10 & Junction & 0\\
\hline
11 & Junction & 0 & 12 & Junction & 0\\
\hline
13 & Junction & 0 & 14 & Junction & 0\\
\hline
15 & Junction & 100 & 16 & Junction & 100\\
\hline
17 & Junction & 100 & 18 & Junction & 0\\
\hline
19 & Junction & 0 & 20 & Junction & 0\\
\hline
21 & Reservoir & 0 & 22 & Reservoir & 0\\
\hline
23 & Tank & 0 & 24 & Tank & 0\\
\hline
25 & Tank & 0 & ~ & ~ & ~\\
\hline
\hline
\end{tabular}
\label{table:nodal_info}
\end{table}

\begin{table*}[htbp!]
\caption{Link Setting of the Barcelona Water Distribution Network (Chezy-Manning)}
\centering
\begin{tabular}{c|c|c|c|c|c|c||c|c|c|c|c|c|c}
\hline
\hline
Link & Type & \tabincell{c}{From\\Node} & \tabincell{c}{To\\Node} & \tabincell{c}{Pipe\\Length\\(feet)} & \tabincell{c}{Pipe\\Diameter\\(feet)} & \tabincell{c}{Pipe\\Rough\\-ness}& Link & Type & \tabincell{c}{From\\Node} & \tabincell{c}{To\\Node} & \tabincell{c}{Pipe\\Length\\(feet)} & \tabincell{c}{Pipe\\Diamter\\(feet)} & \tabincell{c}{Pipe\\Rough\\-ness}\\
\hline
1 & Pipe & 2 & 3 & 2000 & 12 & 0.03 & 2 & Pipe & 3 & 4 & 1000 & 12 & 0.03\\ 
\hline
3 & Pipe & 12 & 15 & 3000 & 12 & 0.03 & 4 & Pipe & 13 & 16 & 4000 & 12 & 0.03\\ 
\hline
5 & Pipe & 14 & 17 & 5000 & 12 & 0.03 & 6 & Pipe & 11 & 18 & 1000 & 12 & 0.03\\ 
\hline
7 & Pipe & 18 & 5 & 1000 & 12 & 0.03 & 8 & Pipe & 23 & 18 & 1000 & 12 & 0.03\\ 
\hline
9 & Pipe & 9 & 19 & 1000 & 12 & 0.03 & 10 & Pipe & 19 & 10 & 1000 & 12 & 0.03\\ 
\hline
11 & Pipe & 24 & 19 & 1000 & 12 & 0.03 & 12 & Pipe & 7 & 20 & 1000 & 12 & 0.03\\ 
\hline
13 & Pipe & 20 & 3 & 1000 & 12 & 0.03 & 14 & Pipe & 20 & 8 & 3000 & 12 & 0.03\\ 
\hline
15 & Pipe & 25 & 20 & 1000 & 6 & 0.03 & 16 & Pump & 21 & 1 & - & - & -\\ 
\hline
17 & Pump & 22 & 6 & - & - & - & 18 & PRV & 6 & 7 & - & - & -\\ 
\hline
19 & PRV & 1 & 9 & - & - & - & 20 & PRV & 1 & 2 & - & - & -\\ 
\hline
21 & PRV & 1 & 11 & - & - & - & 22 & Pipe & 4 & 12 & 3000 & 12 & 0.03\\ 
\hline
23 & Pipe & 5 & 13 & 3000 & 12 & 0.03 & 24 & Pipe & 10 & 14 & 3000 & 12 & 0.03\\ 
\hline
\hline
\end{tabular}
\label{table:link_info}
\end{table*}



\begin{figure}
    \centering
    \subfigure[Water Demand Pattern]{\label{fig:demand_pattern}
    \includegraphics[width=1.6in]{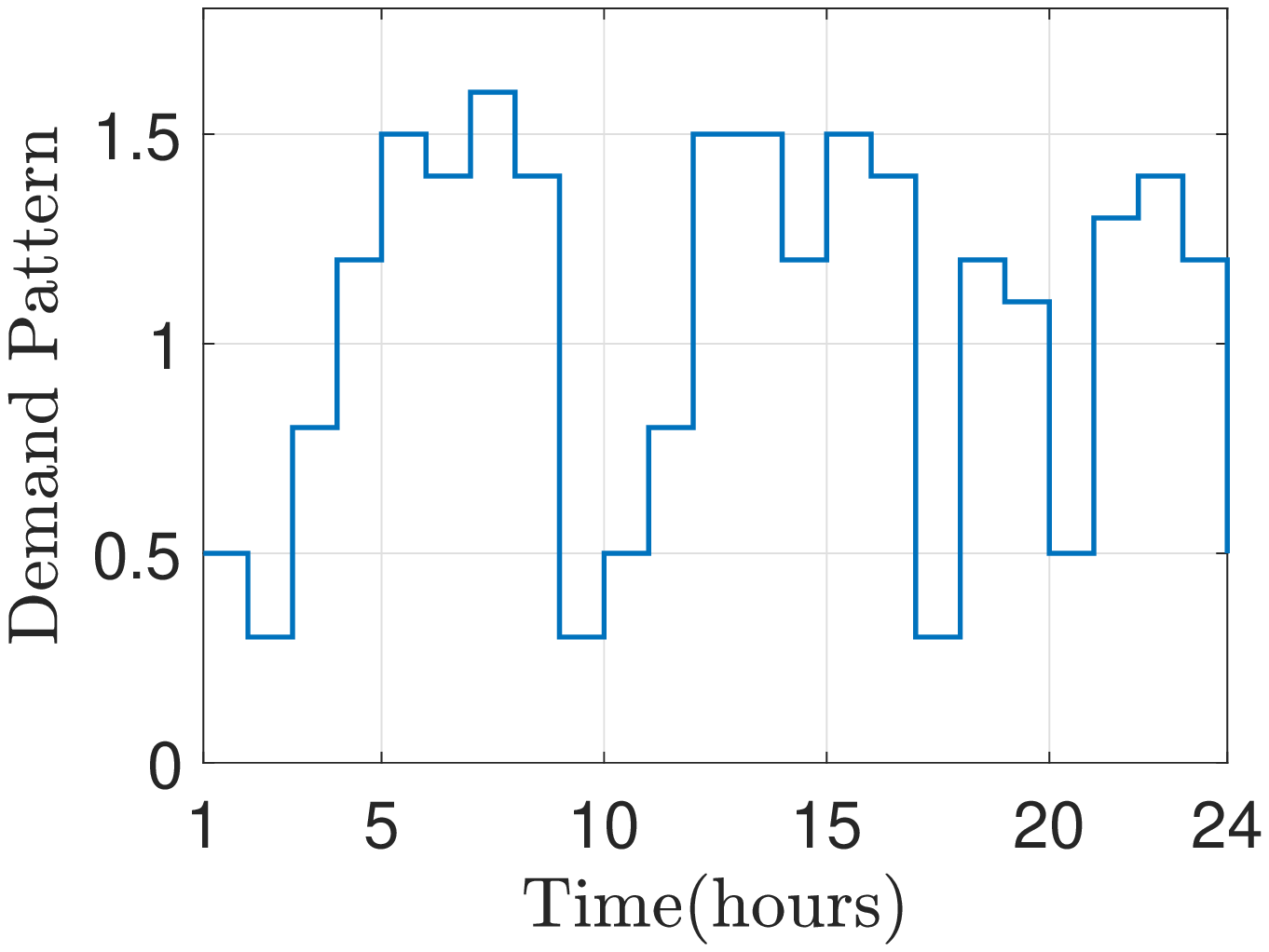}}
    \subfigure[Time-of-Use Price]{\label{fig:TOU}     \includegraphics[width=1.6in]{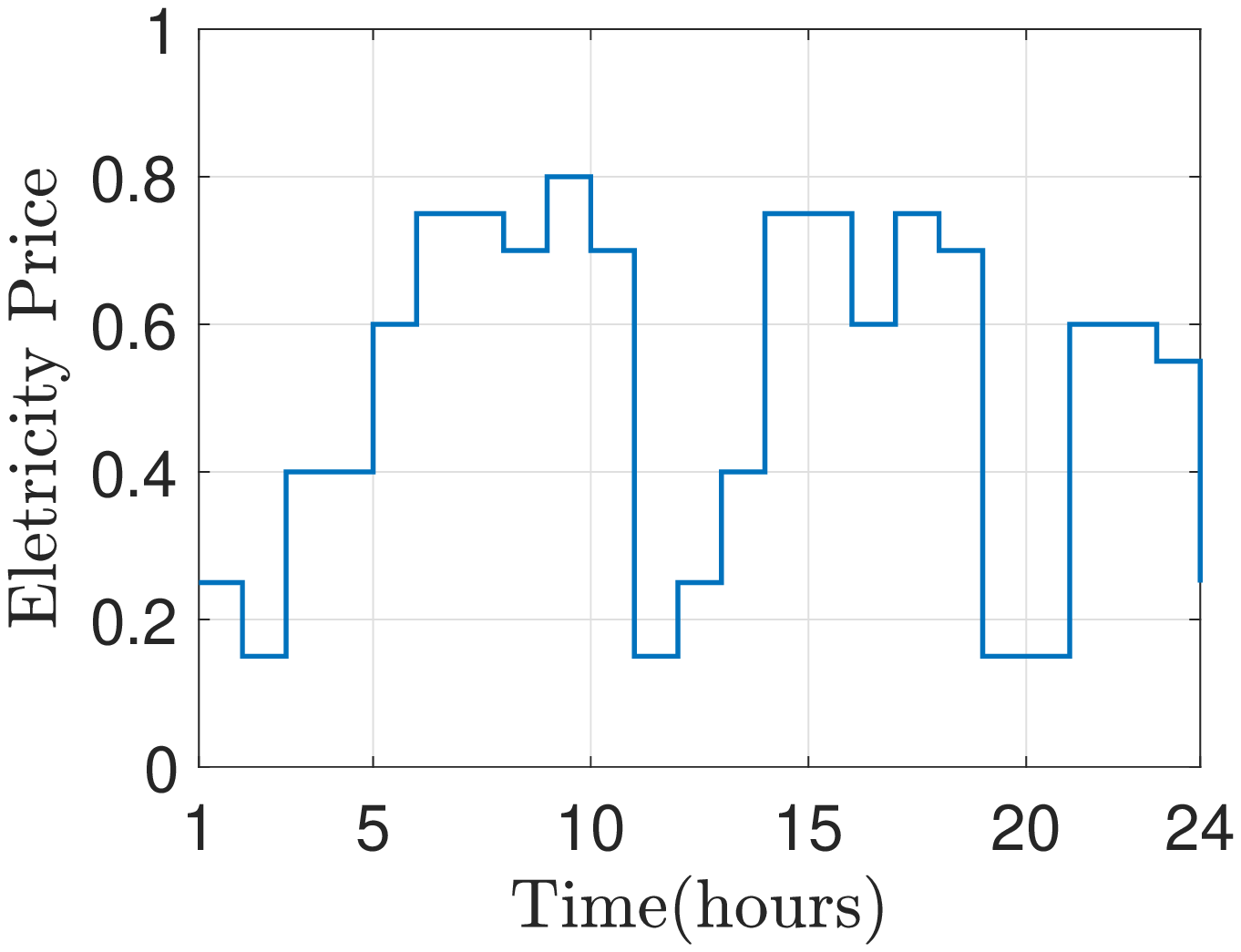}}
    \caption{Time-varying input profiles including water demand pattern and TOU electricity price. The actual water consumption at each node depends on the based demand setting shown in Table \ref{table:nodal_info}. The pump cost function is parameterized in proportional to TOU electricity price.}
\end{figure}

\begin{figure}[!tbhp]
    \centering
    \includegraphics[width=3.5in]{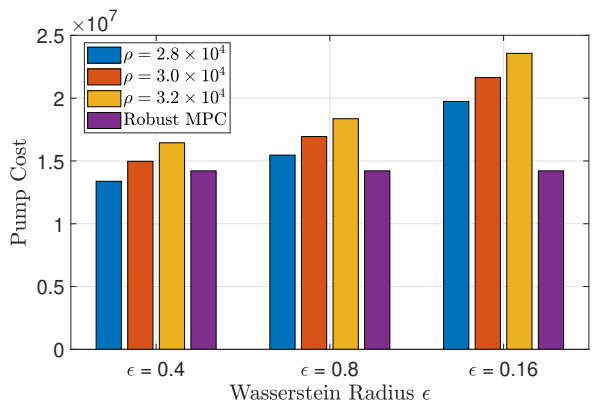}
    \caption{Trade-offs between conservativeness of optimal decisions
    and pump operational costs under various Wasserstein radius $\epsilon$ and risk aversion $\rho$.}
    \label{fig:conservativeness_tradeoff}
\end{figure}

\begin{figure}[!tbhp]
    \centering
    \vspace{2mm}
    \includegraphics[width=3.5in]{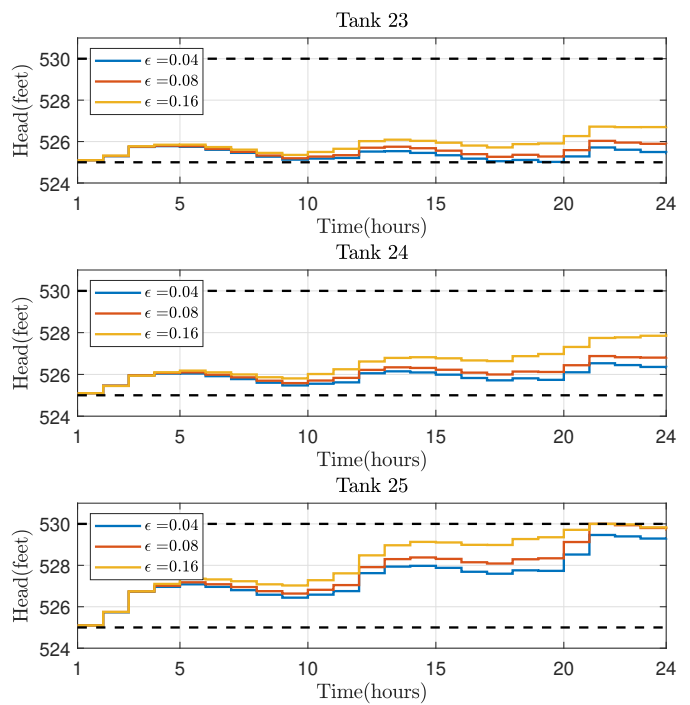}
    \caption{Optimal state trajectories of three tanks (i.e., Node 23, Node 24 and Node 25) for varying Wasserstein radii $\epsilon$. The dash lines indicate the upper and lower bounds on tank head. The initial tank level for all three tanks is 525.1 feet. The risk aversion is set to $\rho = 2.8 \times 10^4$.}
    \label{fig:tank_level}
\end{figure}

\begin{figure}[!tbhp]
    \centering
    \includegraphics[width=3.5in]{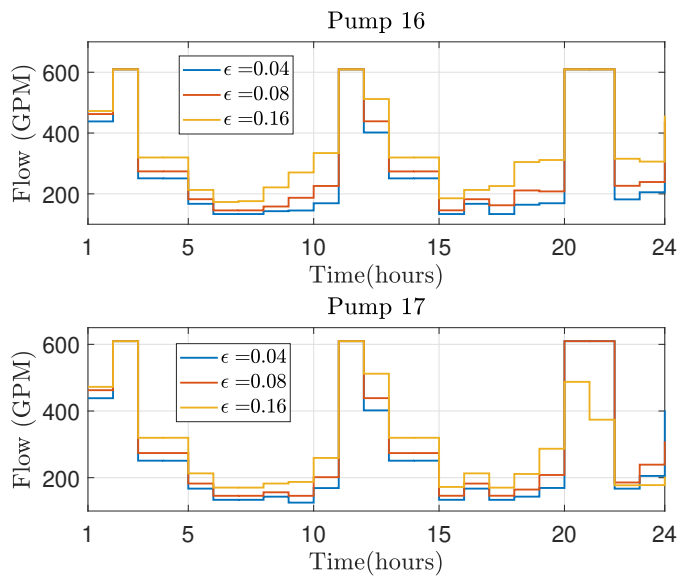}
    \caption{Comparison of optimal pump schedule for various value of Wasserstein radius $\epsilon = 0.04,0.08,0.16$ under certain risk aversion $\rho = 2.8 \times 10^4$.}
    \label{fig:pump_schedule}
\end{figure}

Fig.~\ref{fig:conservativeness_tradeoff} visualizes the fundamental trade-off between the conservativeness of constraint violation and the water network operational costs during 24 hour operation under various risk aversion $\rho$ and Wasserstein radius $\epsilon$. As we increase the Wasserstein radius $\epsilon$, the pump cost will increase as well, but leads to more conservative pumps schedule and lower risk of tank constraint violation. A larger $\epsilon$ results in less constraint violation based on the inherent sampling errors and in turn guarantee the a stronger robustness performance, which will ensure good out-of-sample performance. In addition, with increasing risk aversion $\rho$, the CVaR of constraint violation is emphasized, which comes to a higher operational costs and lower constraint violation. A robust MPC case study is presented by assuming the water demand forecast errors with a support of $\xi \in \Theta := [ -\theta \mathbf{\bar{d}} \leq \xi \leq \theta \mathbf{\bar{d}}, \theta = 0.2]$, which is similar to the Robust MPC formulation in the previous work \cite{goryashko2014robust}. In this comparison study, the linearization nodal pressure constraints were modelled based on the above robust uncertain set and utilized the successive linearization for optimal decisions. Comparison costs are demonstrated in Fig. 4. Data-based DRO OWF provides a controllable conservativeness compared to the stochastic robust MPC. 

Fig.~\ref{fig:tank_level} and Fig.~\ref{fig:pump_schedule} show the water level of tanks $\mathbf{h}_T^\mathrm{TK}$ and the optimal nominal pump schedule $e_T$ over $T - 24$ hours under varying Wasserstein radius $\epsilon$. The tank head trajectories and pump schedule are re-optimized at each time-step via the closed-loop MPC controller based on the data-based uncertainty representation (i.e., a Wasserstein ball of distributions of water demand forecast errors). To prevent the tank level decrease lower then 525 feet, the pumps need to transport more water to tanks for accommodating the water demand uncertainties. As the results, the pumps are more active during the time-slots with higher electricity costs, which leads a significant increase of operational costs. This leads to a safer tank level profiles, as shown in Fig.~\ref{fig:tank_level}. When $\epsilon$ is small, the water consumption mostly come from tanks to maintain an economic operation, which cause the possible constraint violation (e.g., Tank 23 when $\epsilon = 0.04$) if the demand variations were large. As we increase $\epsilon$ leading to a more conservative decision, all pumps sacrifice the operational efficiency and provide more water to increase the tank level and support the water demands. The tank lower level constraints are satisfied due to the better robustness to water demand forecast errors.

To demonstrate the effectiveness of the proposed framework \eqref{eq:DRO_OWF_MPC}, we also introduce the EPANET built-in traditional Rule-based Control (RBC) scheme, which has been widely employed for various water engineering problems. The RBC scheme shares the same control constraints in \eqref{eq:DRO_OWF_MPC}, limits the water heads of three tanks (i.e., Tanks 23, 24 and 25) within a prescribed safe range (i.e., [525, 530] ft) via binary ON or OFF status of pumps (i.e., Pumps 16 and 17). The time step to control pumps is set to one hour, which implies the pumps check the water levels of tanks every hour and then perform control actions. Note that Tanks 23 and 24 can only be controlled by Pump 16 while Tank 25 can be managed by Pumps 16 and 17 simultaneously. Note that the water levels of Tank 23 and Tank 25 possibly direct Pump 16 to the completely conflict control actions (i.e., ON or OFF) if we do not explicitly prioritize these two tanks. Therefore, we assign the level signal from Tank 23 is the priority for Pump 16 to take control actions if a conflict happens.

\begin{figure}[!tbhp]
    \centering
    \vspace{1.5mm}
    \includegraphics[width=3.5in]{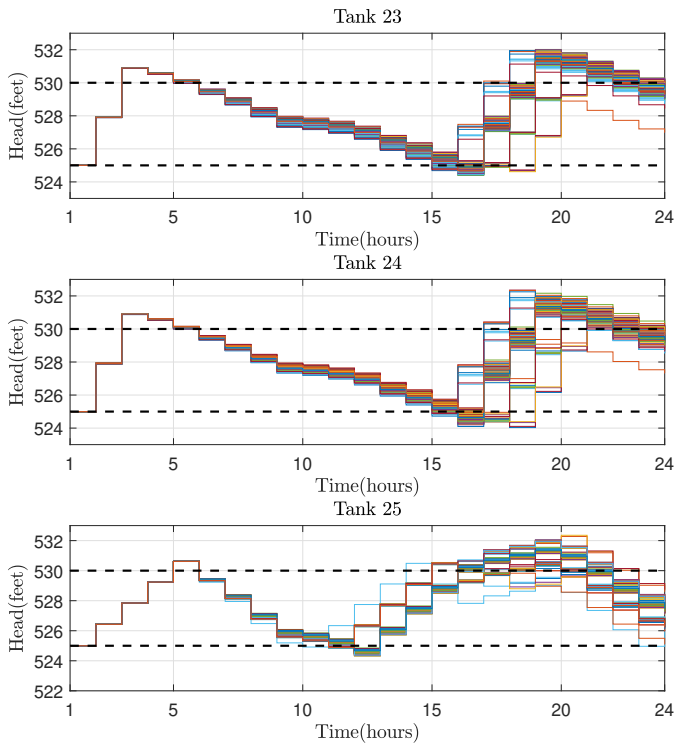}
    \caption{State trajectories of tank head (i.e., Tank 23, Tank 24 and Tank 25) after performing RBC via Monte Carlo simulation including 100 demand scenarios. The dash lines indicate the upper and lower bounds on tank head.}
    \label{fig:tank_level_MC}
\end{figure}

Fig.~\ref{fig:tank_level_MC} illustrates the water levels of three tanks based on the RBC scheme via Monte Carlo simulations. We randomly generate 100 scenarios of water demand forecast errors, which follow the Gaussian distribution with zero mean and 20\% standard deviation of nominal water demand shown in Fig. \ref{fig:demand_pattern}. It readily seen that the RBC mechanism fails to realize the water demand forecast errors and can not successfully manage the tank heads located at a prescribed safe bound. In general, it is very hard to parameterize the RBC control scheme for low risk constraint violation guarantee under the large demand variation, which is due to its decentralized control structure. The benefit of closed-loop multi-period distributionally robust optimal water flow based on model predictive control scheme can be clearly seen via the comparison to the RBC control framework. We analyze the trade-off and conservativeness of stochastic OWF problems in our preliminary works \cite{guo2020distributionally}, including robust MPC, the sample average approximation with CVaR, the chance-constraints (moment-based distributionally robust optimization) MPC and the deterministic scenario. A useful insight that comes out of our results is that the trajectories of tank heads follow similar patterns under varying MPC formulations but with different conservativeness.

In summary, we conclude that our proposed data-based distributionally robust OWF framework can explicitly incorporate water demand uncertainties and successfully control the trade-off between operational efficiency, risk of constraint violation and out-of-sample performance. 

\section{Conclusion and Outlook}\label{sec:conclusion}
We propose a data-based distributionally stochastic robust optimal water flow based on limited information from water demand forecasts. The framework creates and then leverages a successive linearization of hydraulic coupling for an efficient computation of multi-period feedback control policies, which are robust to inherent sampling errors in the training dataset. We explore the tractability of proposed multi-period OWF problem via the Wasserstein-based distributional information of ambiguity set centered at the empirical distribution. The effectiveness and flexibility has been demonstrated on a 25-node water distribution network for the optimal water pump schedule and tank head management under water demand uncertainties. The numerical results indicate that our proposed framework has superior out-of-sample performance then existing control frameworks and allows flexible parameterization to systematically exploit the operating strategies of water pumps to explicitly trade-off the operational efficiency and constraint violations due to large water demand variations. The limitations and outlooks of the proposed data-based optimal pump control are summarized below:

\underline{Scalability:} As the size of the water network and the number of water demand samples increase, more computational efforts are required to solve the optimal solutions. Though the proposed optimization problem is convex and can be solved by many solvers, future work will focus on improving the computational affordability of our framework for large-scale water distribution networks \cite{muros2018game}.
            
\underline{Communication:} The optimal pump control actions are attained by successively solving a distributionally robust optimization, which requires global communication. Future work will extend the proposed framework to a distributed setting with local communication.
            
\underline{Network Topology:} The Laplacian approximation for solving intractable optimal water flow works appropriately for water distribution networks with pure tree topologies. A more general successive Laplacian approximation is required for mesh water distribution networks in the future.
            
\underline{Feasibility:} The feasibility of the proposed data-based distributionally robust stochastic MPC is numerically validated instead of theoretically. Further feasibility analysis of a generic data-based distributionally robust stochastic MPC and the convergence of the successive algorithm will be theoretically established.\\
Future works also can be
\begin{itemize}
    \item developing a data-based distributionally robust control framework for optimal water contamination control;
    \item including operational status of actuators as controllable variables for a distributionally robust stochastic hybrid MPC OWF framework;
    \item developing a theoretical convergence analysis of the proposed successive linearization algorithm.
\end{itemize}

\section*{Appendix A}
\begin{figure}[h!]
        \centering
        \includegraphics[width=3.4in]{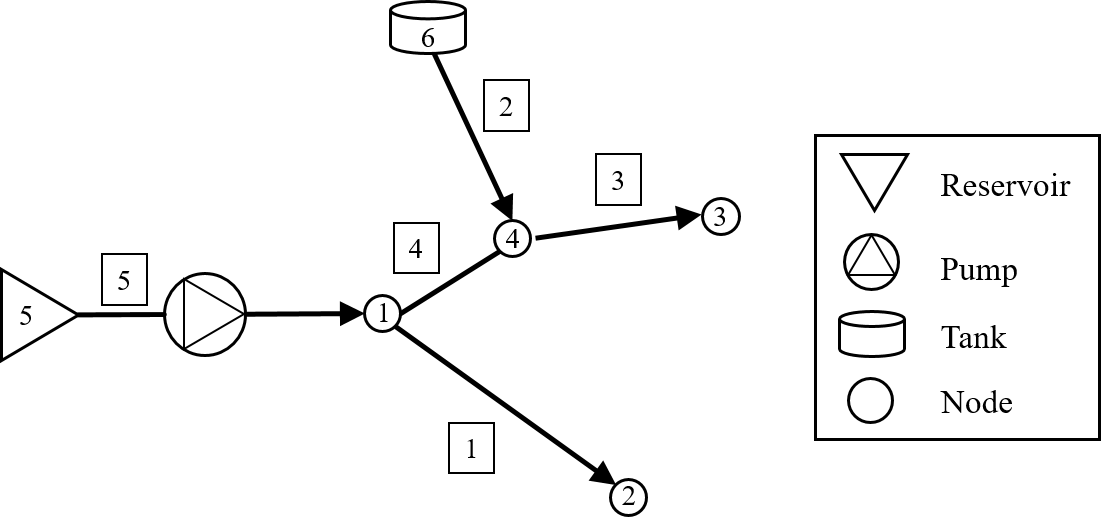}
    \caption{A 6-node water distribution network with 1 reservoir, 1 tank, 1 pump, and 2 water demands.}
    \label{fig:tree6}
    \end{figure}

Here, we demonstrate a comparison between the feasibility results from the proposed Laplacian-based successive linearization approach in Algorithm \ref{algorithm:successive} and the water flow results from the EPANET nonlinear solver. For simple demonstration, we used a 6-node water distribution network with 1 reservoir, 1 tank, 1 pump and 2 demands in Fig.~\ref{fig:tree6}. Two demands are located in Node 2 and Node 3 with the base demand of 50 GPM. The nominal water demand pattern over 24 hours is the same as shown in Fig.~\ref{fig:demand_pattern}. 
\color{black}

\begin{figure}[ht!]
        \centering
        \includegraphics[width=3.5in]{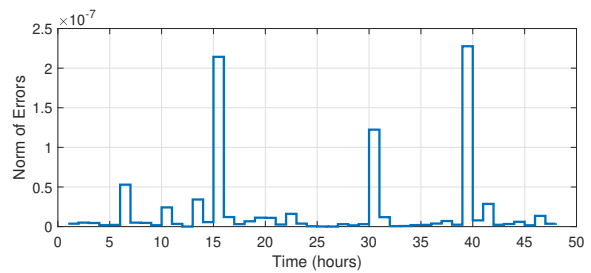}
    \caption{The linearized errors of the successive algorithm over 24 hours, defined by $\textrm{Err}:=\|\Delta \mathbf{q}_n^*\|_2^2$.}
    \label{fig:linear_errors}
    \end{figure}

\begin{figure*}[ht!]
        \centering
        \includegraphics[width=7.5in]{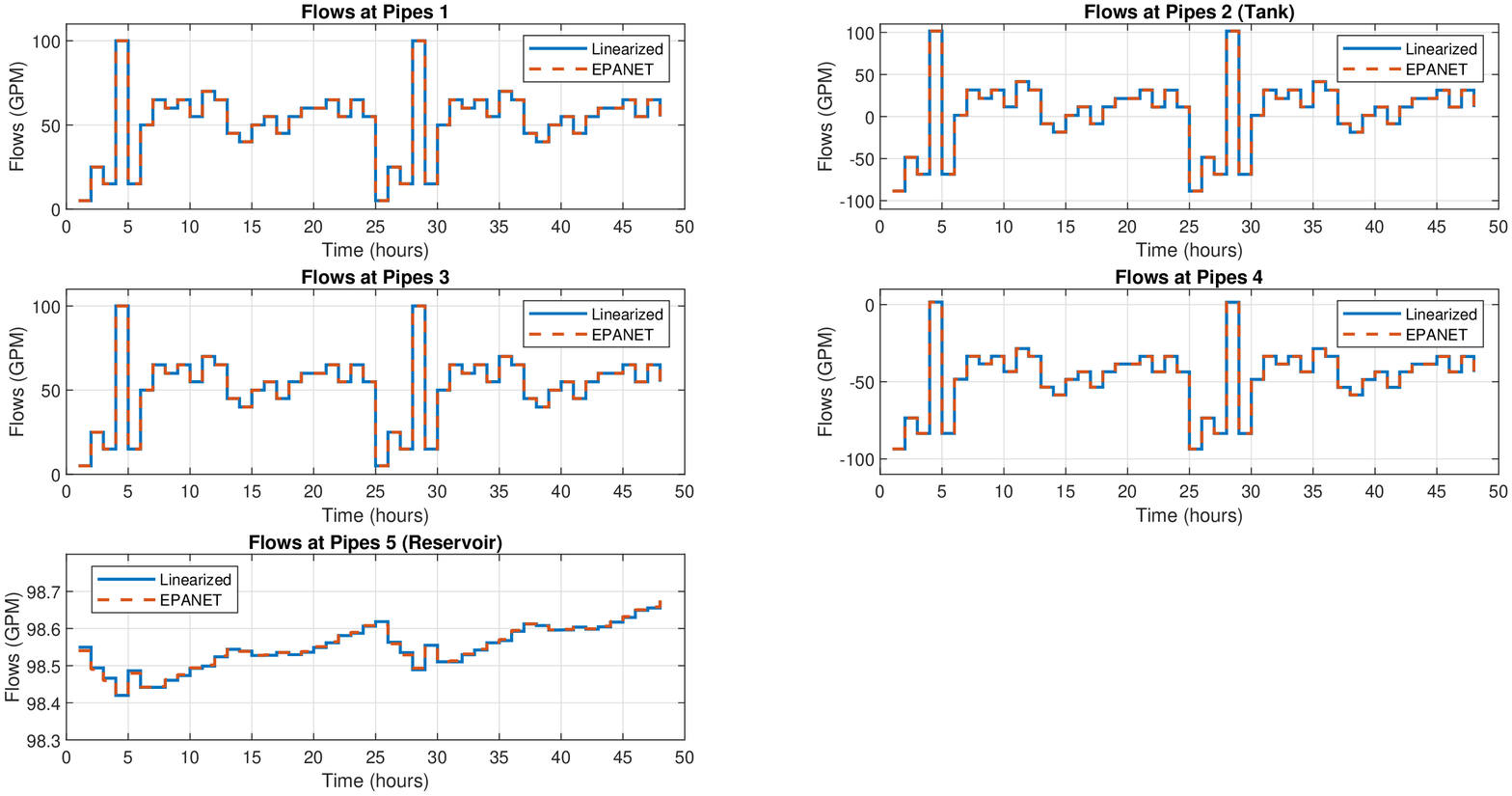}
    \caption{Comparison of water flow status between the feasibility results from the proposed successive Laplacian-based linearized approach and the results from the EPANET nonlinear solver.}
    \label{fig:linear_Q}
    \end{figure*}

\begin{figure*}[ht!]
        \centering
        \includegraphics[width=7.5in]{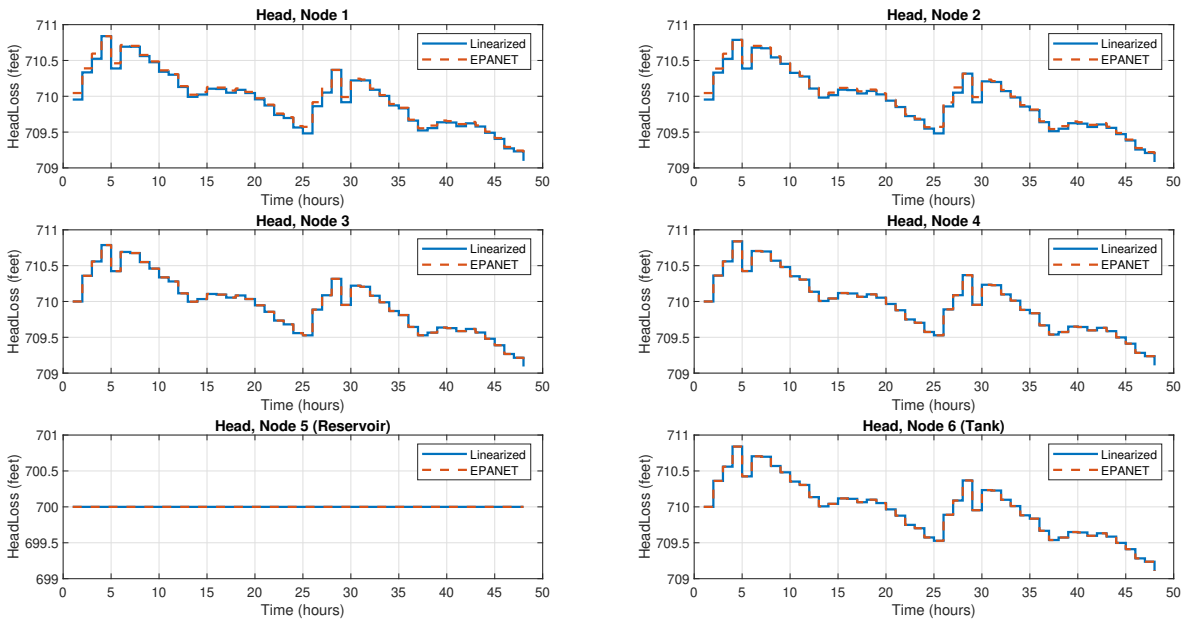}
    \caption{Comparison of head status between the feasibility results from the proposed successive Laplacian-based linearized approach and the results from the EPANET nonlinear solver.}
    \label{fig:linear_H}
    \end{figure*}

Fig.~\ref{fig:linear_Q} and Fig.~\ref{fig:linear_H} visualize the results from the comparison, which demonstrate the effectiveness of the proposed Laplacian-based successive linearization approximation. Fig.~\ref{fig:linear_errors} also presents the linearization errors by implementing Algorithm \ref{algorithm:successive}, which is defined by $\textrm{Err}:=\|\Delta \mathbf{q}_n^*\|_2^2$. The maximum linearized error over 24 hours is less than $2.5\times 10^{-7}$. It can be seen that the water head and flow trajectories from the Laplacian-based model successfully provide a great approximation of the network status solved by the exact nonlinear hydraulic models.

\section*{Appendix B}
Here we provide the constant matrices $\{A,B,E,F\}$ shown in \eqref{eq:State_space_expression} in terms of the network coefficients in the Laplacian-based linearization \eqref{eq:Laplacian_Energy_Conservation}.
Recall the tank level dynamics \eqref{eq:tank_SS} 
\begin{equation*}
x(t+1) = \bar{A} x(t) + \bar{B}_uu(t) + \bar{B}_v v(t),
\end{equation*}
where $\bar{A}$ is the $n_{\textrm{TK}}$-dimension identity matrix and the last two terms on the right side are re-organized in terms of the network incident matrix $\mathbf{B}_f^\intercal$, the coefficient of tanks $A_i^{\mathrm{TK}}$ and time interval $\Delta t$
\begin{equation}\nonumber
    \begin{bmatrix}
    \bar{B}_u & \bar{B}_v
    \end{bmatrix}
    \begin{bmatrix}
    u(t)\\
    v(t)\\
    \end{bmatrix}= \underbrace{
    \begin{bmatrix}
    \frac{\Delta t}{A_1^{\textrm{TK}}} & \cdots  & 0\\
    \vdots & \ddots & \vdots\\
    0 & \cdots & \frac{\Delta t}{A_{n_\textrm{TK}}^{\textrm{TK}}}
    \end{bmatrix}[\mathbf{B}_f^\intercal]_{\{\mathcal{T}\}}}_{\begin{bmatrix}
    \bar{B}_u & \bar{B}_v
    \end{bmatrix}}
    \begin{bmatrix}
    u(t)\\
    v(t)\\
    \end{bmatrix}.
\end{equation}
The operator $[\cdot]_{\{\mathcal{T}\}}$ selects the $i^\mathrm{th}$-row of matrix $\mathbf{B}_f^\intercal$, where all node $i$ are collected in the set $\mathcal{T}$. Similarly, for mass balance
    $d(t) = \bar{E}_u u(t) + \bar{E}_v v(t)$
we have
\begin{equation}\nonumber
    \begin{bmatrix}
    \bar{E}_u & \bar{E}_v
    \end{bmatrix}
    \begin{bmatrix}
    u(t)\\
    v(t)
    \end{bmatrix} = \mathbf{B}_f^\intercal \begin{bmatrix}
    u(t)\\
    v(t)
    \end{bmatrix}.
\end{equation}
Finally, recall the mass conservation in DAE model \eqref{eq:energy_conservion_SS}
\begin{equation*}
\bar{F}_x x(t) + \bar{F}_l l(t) = \bar{F}_u u(t) + \bar{F}_v v^{\mathrm{P}}(t) + \bar{F}_\phi \phi(t) + F_0,
\end{equation*}
and re-write the above equation as follows
\begin{equation*}
    \underbrace{\begin{bmatrix}
    \bar{F}_x & \bar{F}_l
    \end{bmatrix}}_{\mathbf{L}}
    \underbrace{\begin{bmatrix}
    x(t)\\
    l(t)
    \end{bmatrix}}_{\mathbf{h}(t)}
= 
\underbrace{
\begin{bmatrix}
\bar{F}_u & \bar{F}_v & \bar{F}_{\phi}
\end{bmatrix}}_{\mathbf{B}_f^\intercal}
\underbrace{
\begin{bmatrix}
u(t)\\
v^{\mathrm{P}}(t)\\
\phi(t)
\end{bmatrix}}_{\Delta\mathbf{q}(t)}
+ \underbrace{F_0}_{\bar{\mathbf{Q}}(t)}.
\end{equation*}

\bibliographystyle{ieeetr}        
\bibliography{reference}           

\end{document}